\documentclass[12pt]{article}
\usepackage{amsmath,amssymb}
\usepackage{bm}
\usepackage[usenames]{color}
\newtheorem{tm}{Theorem}[section]
\newtheorem{lm}[tm]{Lemma}

\newtheorem{re}[tm]{Remark}

\newtheorem{pr}[tm]{Proposition}
 
 \newenvironment{demo}[1]{\par\smallskip\par\begin{trivlist}
\item[]{\bf #1}\ }{\end{trivlist}\par\smallskip\par}
\newcommand{\Proof}{\begin{demo}{{\it Proof.\ }}}
\newcommand{\qed}{\end{demo}}
\newcommand{\toy}{\ \rule[0em]{0.5ex}{1.8ex}}
\newcommand{\QED}{\toy\end{demo}}

\newcommand{\la}{\langle}
\newcommand{\ra}{\rangle}
\newcommand{\nn}{\nonumber}
\newcommand{\III}{{\vert \kern-.10em \vert \kern-.10em \vert}}
\newcommand{\ve}{\varepsilon}
\newcommand{\ld}{\lambda}
\newcommand{\al}{\alpha}

\makeatletter
 
 \@addtoreset{equation}{section}
\makeatother
 
\begin{document}
\setlength{\baselineskip}{15pt} 
%
\bibliographystyle{plain}
\title{
Short time kernel asymptotics for rough differential equation
driven by fractional Brownian motion
\footnote{
{\bf Mathematics Subject Classification}: 60H07, 60F99, 60G22.
{\bf Keywords}: rough path theory, Malliavin calculus, 
 fractional Brownian motion, short time asymptotic expansion}
}
\author{
Yuzuru INAHAMA
\footnote{
Graduate School of Mathematics,   Kyushu University
Motooka 744, Nishi-ku, Fukuoka 819-0035, Japan.
E-mail:~\tt{inahama@math.kyushu-u.ac.jp}
}
}
\date{   }
%
%
\maketitle

%
%
\begin{abstract}
We study a stochastic differential equation
in the sense of rough path theory 
driven by fractional Brownian rough path with Hurst parameter $H ~(1/3 < H \le 1/2)$
under the ellipticity assumption at the starting point.
In such a case, the law of the solution at a fixed time has a kernel, i.e., 
a density function with respect to Lebesgue measure.
In this paper we prove a short time off-diagonal asymptotic expansion of the kernel 
under mild additional assumptions.
Our main tool is Watanabe's distributional Malliavin calculus.
\end{abstract}

%
%
\section{Introduction }

For the usual $d$-dimensional Brownian motion $(w_t)$ and 
sufficiently regular vector fields $V_i~(0 \le i \le d)$ on ${\mathbb R}^n$,
consider the following stochastic differential equation (SDE)
of Stratonovich type:
\begin{equation}
dy_t =  \sum_{i=1}^d  V_i (y_t) \circ dw^i_t + V_0 (y_t) dt \qquad  \mbox{ with }  
\qquad y_0 =a \in {\mathbb R}^n. 
\nn
\end{equation}
If the vector fields satisfy the hypoellipticity condition 
at the starting point $a$,
then the law of $y_t$ has a heat kernel i.e., a density function $p_t (a,a')$
with respect to Lebesgue measure $da'$ for any $t>0$.

In probability theory, the short time asymptotic (off-diagonal)
problem of $p_t (a,a')$ has extensively been studied
and is now a classical topic.
See for instance \cite{az, bena1, bena2, bena3, bl1, bl2, bi, gav, ks1, ks2, 
lean1, lean2, lean3, lean4, lean5, mol, tak, tw, ue1, ue2, uw, wa} 
and references therein.
(There are also analytic approaches, of course. But, we do not discuss them in this paper.)
Among many probabilistic methods, Malliavin calculus is known to be quite powerful. 
Bismut \cite{bi} was first to prove short time kernel asymptotics via
 Malliavin calculus.
Among such proofs, we focus on Watanabe's theory of generalized Wiener functionals 
and asymptotic theorems for them \cite{wa, iwbk, tw}.

Recently, the theory of "SDE" for fractional Brownian motion (fBm) was developed.
As a result,  an analogous asymptotic problem is gathering attention.
When Hurst parameter $H$ is larger than $1/2$,
the SDE above is in the sense of Young integration.
When $1/4 < H \le 1/2$,
it should be understood as a differential equation in the rough path sense 
driven by  fractional Brownian rough path. 
In his previous paper \cite{ina1}, the author studied both on-diagonal and off-diagonal
short time asymptotic expansion of $p_t (a, a')$ when $H >1/2$.
The method is Watanabe's asymptotic theory of generalized Wiener functionals
(i.e., Watanabe distributions) in \cite{wa}.
In \cite{ina1} the coefficient vector fields are
assumed to satisfy the ellipticity condition  at $a$ 
and some additional mild conditions are also assumed.
Those conditions are almost parallel to the ones in \cite{wa}.
Simply put, \cite{ina1} is  a "fractional version"
of \cite{wa} in the framework of Young integration.

The aim of this paper is to prove a similar off-diagonal asymptotic expansion 
when $1/3 < H \le 1/2$.
Although the basic strategy of proof is similar to the case $H>1/2$ in \cite{ina1},
the proof gets much more technically difficult 
since we work on the rough path space.
We will carry it out
 by combining various recently proven results for Gaussian rough paths.
A number of paper have been published 
on Malliavin calculus for Gaussian rough paths by now.
See \cite{ai, bnot, boz1, boz2, bgq, cf, cfv, chlt, dri, hp, ht, ina06, ina2, ina3} for instance.
However, this type of short time kernel asymptotics seems new.

The organization of this paper is as follows:
In Section \ref{sec.result}
we give a precise formulation of our problem and
the statement of our main result (Theorem \ref{thm.MAIN.off}).
In Section \ref{sec_two} we prove moment estimates for 
Taylor expansion of Lyons-It\^o map.
The expansion in the deterministic sense is already known,
but we need "$L^p$-version" (or "${\bf D}_{\infty}$-version") of the expansion in this paper.
These estimates play a crucial role in the proof of the main theorem.
In Section \ref{sec.malli} 
we present
two propositions (Propositions \ref{pr.d_infty} and \ref{pr.nondeg}) 
on regularity 
in the sense of Malliavin calculus of the solution of RDE
driven by fractional Brownian rough path.
Thanks to these propositions,  we can use  Watanabe's asymptotic theory
in the proof of the main theorem
in Section \ref{sec.pf}, following the argument in \cite{wa, ina1}.
A difference from \cite{wa} is that 
we can work
and, in particular, localize around the energy minimizing path
in the domain (not in the range) of Lyons-It\^o map since the map is continuous in rough path theory.
%
%
%
%

We do not give a heuristic sketch of our argument for brevity. 
Since formal computations are basically the same as in the Young case,
the reader who wants to know it may consult
the corresponding part of the author's previous paper \cite{ina1}.

\begin{re}
The first version of this paper contains detailed proofs of 
Theorem 3.4, Proposition 4.1,  Proposition 4.2, and Lemma 5.2.
It can be found on arXiv preprint server (arXiv:1403.3181). 
\end{re}

\noindent 
{\bf Acknowledgement}~The author thanks Dr. Nobuaki Naganuma 
for helpful comments.




\section{Setting and main results}
\label{sec.result}

\subsection{Setting}
In this subsection, we introduce a stochastic process that will play a main role in this paper.
From now on we denote by
$w =(w_t)_{ t \ge 0} =( w^1_t, \ldots, w^d_t)_{  t \ge 0}$
 the $d$-dimensional fractional Brownian motion (fBm) with Hurst parameter $H$.
Throughout this paper we assume $1/3 <H \le 1/2$.
It is a unique $d$-dimensional, mean-zero, continuous Gaussian process with covariance
$$
{\mathbb E} [  w^i_s  w^j_t] = \frac{  \delta_{ij} }{2} ( |s|^{2H} +  |t|^{2H} -  |t-s|^{2H}), 
\qquad 
(s, t \ge 0).
$$
Note that, for any $c>0$, $(w_{ct})_{ t \ge 0}$ and $(c^{H} w_t)_{ t \ge 0}$ have the same law.
This property is called self-similarity or scale invariance.
When $H =1/2$, it is the usual Brownian motion.
It is well-known that $w$ admits 
a canonical rough path lift ${\bf w}$, which is called fractional Brownian rough path.

Let $V_i : {\mathbb R}^n \to {\mathbb R}^n ~$ be $C_b^{\infty}$,
that is, $V_i$ is a bounded smooth function with 
bounded derivatives of all order ($0 \le i \le d$).
We consider the following rough differential equation (RDE);
\begin{equation}\label{main.ygODE.eq}
dy_t =  \sum_{i=1}^d  V_i (y_t) dw^i_t + V_0 (y_t) dt \qquad  \mbox{ with }  
\qquad y_0 =a \in {\mathbb R}^n. 
\end{equation}
This RDE is driven by the Young pairing 
$({\bf w}, \bm{\lambda})$, where $\lambda_t =t$. 
The unique solution is denoted by ${\bf y} =({\bf y}^1, {\bf y}^2)$
and we set $y_t := a + {\bf y}_{0,t}^1$ as usual.
We will sometimes write $y_t =y_t(a) = y_t(a, {\bf w})$ etc. 
to make explicit the dependence on $a$ and ${\bf w}$.

A matrix notation is often convenient.
So we set $b =V_0$
and $\sigma =[V_1, \ldots, V_d]$, which is $n \times d$ matrix-valued,
and often rewrite RDE (\ref{main.ygODE.eq}) as follows;
\[
dy_t =  \sigma (y_t) dw_t + b (y_t) dt \qquad  \mbox{ with }  
\qquad y_0 =a \in {\mathbb R}^n.\]

\subsection{Assumptions}

In this subsection we introduce assumptions of the main theorems.
First, we assume the ellipticity of the coefficients of (\ref{main.ygODE.eq}) at the starting point $a \in {\mathbb R}^n$.
\\
\\
{\bf (A1):}~ 
The set of vectors $\{V_1(a), \ldots, V_d(a)\}$ linearly spans ${\mathbb R}^n$.
\\
\\
It is known that, under Assumption {\bf (A1)}, 
the law of the solution $y_t$ has a density $p_t(a, a')$ with respect to the Lebesgue measure 
on ${\mathbb R}^n$
for any $t >0$ (see \cite{hp}).
Hence, for any Borel subset $U \subset {\mathbb R}^n$, 
${\mathbb P} (y_t (a) \in U) = \int_U p_t(a, a') da'$.

Let ${\cal H} = {\cal H}^H$ be the Cameron-Martin space of fBm $(w_t)$.
Note that 
any $\gamma \in {\cal H}$ is continuous and of finite $q$-variation 
for some $q \in [1,2)$.
For $\gamma \in {\cal H}$, we denote by $\phi^0_t  =\phi^0_t (\gamma)$ be the solution
of the following Young ODE;
\begin{equation}\label{main.phi0.eq}
d\phi^0_t 
=  \sum_{i=1}^d  V_i (  \phi^0_t) d\gamma^i_t
 \qquad  \mbox{ with }  \qquad \phi^0_0 =a \in {\mathbb R}^n. 
\end{equation}
Set,  for $a' \neq a$, 
\[
K_a^{a'} = \{ \gamma \in {\cal H} ~|~  \phi^0_1(\gamma) =a'\}.
\]
We only consider the case where $K_a^{a'} $ is not empty.
For example, if we assume {\bf (A1)} for all $a$, then this set $K_a^{a'} $ is not empty.
%
%
%
From  goodness of the rate function in Schilder-type large deviation
for fractional Brownian rough path (see \cite{fv07}), 
it follows that
$\inf\{ \|\gamma\|_{\cal H} ~|~ \gamma \in  K_a^{a'}\} 
= \min\{ \|\gamma\|_{\cal H} ~|~ \gamma \in  K_a^{a'}\}$.
Now we introduce  the following assumption;
\\
\\
{\bf (A2):} $\bar{\gamma} \in K_a^{a'}$ which minimizes ${\cal H}$-norm exists uniquely.
\\
\\
In what follows, $\bar{\gamma}$ denotes the minimizer in Assumption {\bf (A2)}.
We also assume that 
$\| \,\cdot\, \|^2_{{\cal H}} /2$ is not so degenerate at $\bar\gamma$ in the following sense.
\\
\\
{\bf (A3):} At $\bar{\gamma}$, the Hessian of the functional  
$K_a^{a'} \ni \gamma \mapsto \| \gamma \|^2_{{\cal H}} /2 $ is strictly positive in the quadratic form sense.
More precisely, 
if $(- \ve_0, \ve_0) \ni u \mapsto f(u) \in K_a^{a'}$ is a smooth curve in $K_a^{a'}$
such that $f(0) = \bar\gamma$ and $f^{\prime} (0)  \neq 0$, then
$
(d/du)^2 \vert_{u=0} \|   f(u) \|^2_{{\cal H}} /2   > 0.$
\\
\\
Later we will give a more analytical condition {\bf (A3)'}, which is equivalent to {\bf (A3)} 
under {\bf (A2)}.
In \cite{wa}, Watanabe used {\bf (A3)'} in his proof of off-diagonal kernel asymptotics.
We will also use {\bf (A3)'}.
In order to state {\bf (A3)'}, however, we have to introduce a lot of notations.
So, we presented {\bf (A3)} here for ease of presentation.

\begin{re}
Assume {\bf (A1)}.
If the end point $a'$ is sufficiently close to the starting point $a$, 
then {\bf (A2)} and {\bf (A3)} are satisfied.
(This is shown in the author's previous paper \cite{ina1}
when $1/2 <H<1$.
The same proof  works in our case ($1/3 <H\le 1/2)$, too.
The key is the implicit function theorem.)
\end{re}

\subsection{Index sets}

In this  subsection we introduce several index sets for the exponent of the small parameter $\ve >0$, 
which will be used in the asymptotic expansion.
Unfortunately,  index sets in this paper are not 
the set of (a constant multiple of) natural numbers 
and are rather complicated.
(However, all these index sets are discrete subsets of $({\mathbb Z}+ H^{-1}{\mathbb Z}) \cap [0,\infty)$
with the minimum $0$.)

Set
\[
\Lambda_1 = \{   n_1 + \frac{n_2}{H}  ~|~  n_1, n_2 \in {\mathbb N}  \},
\]
where ${\mathbb N}= \{0,1,2,\ldots\}$.
We denote by $0=\kappa_0 <\kappa_1 < \kappa_2 < \cdots$ all the elements of $\Lambda_1$ in increasing order.
For a while, consider the case $1/3 <H <1/2$.
Several  smallest elements are explicitly given as follows;
\[
\kappa_1 = 1, \quad   \kappa_2 = 2, \quad  
\kappa_3 =\frac{1}{H},  \quad   \kappa_4 = 3, \quad 
\kappa_5 = 1+\frac{1}{H},  
\quad  \kappa_6 = 4, \ldots
\]
As usual,  using the scale invariance (i.e., self-similarity) of fBm,
we will consider the scaled version of (\ref{main.ygODE.eq}).
(See the scaled and shifted RDE (\ref{rde_mal.eq}) below).
From its explicit form, one can easily guess why $\Lambda_1$ appears. 

We also set 
\[
\Lambda_2 = \{ \kappa -1 ~|~ \kappa \in \Lambda_1 \setminus \{0\} \}
=
\Bigl\{ 
0,  \, 1, \,
\frac{1}{H} -1,  \, 2,   \, \frac{1}{H}, 
  \, 3,  \ldots
 \Bigr\}
\]
and 
\[
\Lambda'_2 = \{ \kappa -2 ~|~ \kappa \in \Lambda_1 \setminus \{0, 1\} \}
=
\Bigl\{ 
0,   \, \frac{1}{H} -2, 
  \, 1, \, \frac{1}{H}-1, 
  \, 2,  \ldots
 \Bigr\}.
\]
Next we set 
\[
\Lambda_3 = 
\{  a_1 +a_2 + \cdots + a_m ~|~  \mbox{$m \in {\mathbb N}_+$ and $a_1 ,\ldots, a_m \in \Lambda_2$} \}.
\]
In the sequel, $\{ 0=\nu_0 <\nu_1<\nu_2 <\cdots \}$ stands for all the elements of 
$\Lambda_3$ in increasing order.
Similarly, 
\[
\Lambda'_3 = 
\{  a_1 +a_2 + \cdots + a_m ~|~  \mbox{$m \in {\mathbb N}_+$ and $a_1 ,\ldots, a_m \in \Lambda'_2$} \}.
\]
In the sequel, $\{ 0=\rho_0 <\rho_1<\rho_2 <\cdots \}$ stands for all the elements of 
$\Lambda'_3$ in increasing order.
Finally,  
$$\Lambda_4 =\Lambda_3  +\Lambda'_3 =\{ \nu + \rho ~|~ \nu   \in  \Lambda_3  , \rho \in  \Lambda'_3\}.$$
We denote by $\{ 0 = \lambda_0 < \lambda_1 < \lambda_2 < \cdots \}$ all the elements of $\Lambda_4$  in increasing order.

When $H =1/2$, all these index sets $\Lambda_i, \Lambda'_j$ above are just ${\mathbb N}$.

\subsection{Statement of the main result}

Now we state our main theorem,
which is basically analogous to the corresponding one in Watanabe \cite{wa}.
However, when $H \neq 1/2$ and the drift term exists,
there are some differences.
First,  the exponents of $t$ are not (a constant multiple of) natural numbers.
Second, cancellation of "odd terms" as in p. 20 and p. 34, \cite{wa}
does not occur in general. 
(These phenomena were already observed in \cite{ina1}  in the Young integration setting i.e., the case $H >1/2$.)

\begin{tm} \label{thm.MAIN.off}
Assume $a \neq a'$ and {\bf (A1)}--{\bf (A3)}. 
Then, we have the following asymptotic expansion as $t \searrow 0$;
    \[
p(t , a, a')  \sim
    \exp \Bigl(  -\frac{ \|  \bar\gamma \|^2_{{\cal H}}}{2  t^{2H} }  
     \Bigr)
  \frac{1}{t^{n H}} 
     \bigl\{
      \alpha_{ 0}  +  \alpha_{\lambda_1}  t^{ \lambda_1 H} + \alpha_{\lambda_2}  t^{\lambda_2  H  }  +\cdots
            \bigr\}
\]
 for certain real constants $\alpha_{\lambda_j}  ~ (j =0,1,2,\ldots)$. 
 Here, $\{ 0=\lambda_0 <\lambda_1<\lambda_2 <\cdots \}$ are all the elements of 
$\Lambda_4$ in increasing order.
Moreover, $\alpha_0$ is positive.
 \end{tm}

\begin{re}\label{re.const}
{\bf (i)}~
In theory, the constants in the asymptotic expansion in 
Theorem \ref{thm.MAIN.off} (and in the on-diagonal case in Theorem \ref{tm.ondiag} below) 
are computable.
But, actual computation is quite cumbersome and we do not carry it out in this paper.
We just mention here that the first constants $\alpha_{ 0}$ in  
Theorem \ref{thm.MAIN.off} and $c_{0} $ in Theorem \ref{tm.ondiag} are non-zero.
\\
\noindent
{\bf (ii)}~
It might be interesting to consider the case $1/4<H \le 1/3$.
In that case, since the third level rough path theory is needed,
calculations may become much harder.
\\
\noindent
{\bf (iii)}~Our assumptions {\bf (A1)}--{\bf (A3)} are quite similar to 
the corresponding ones in \cite{wa}.
Therefore, if we set $H=1/2$ in 
 Theorem \ref{thm.MAIN.off} above recovers most 
 of (but not all of) the main result in Watanabe \cite{wa}.
Hence, our result could also be regarded as a rough path proof of \cite{wa}.
(In this case, however, the index set in Theorem \ref{thm.MAIN.off} is not $\Lambda_4 ={\mathbb N}$,
but is actually $2{\mathbb N}$, due to cancellation of the odd terms.)
Compared to the main theorem in \cite{wa},
Theorem \ref{thm.MAIN.off} with $H=1/2$ does not include the following two cases;
\\
{\rm (a):}~In this paper the ellipticity assumption {\bf (A1)} is assumed.
In \cite{wa}, however, something like "step 2-hypoellipticity" case was also studied.
(We simply did not try this case.)
\\
{\rm (b):}~In this paper the coefficient vector fields are of $C_b^{\infty}$.
However, the condition on vector fields in \cite{wa} is as follows:
"For all $m=1,2,\ldots$ and $0 \le i \le d$,
$ \|\nabla^m V_i \|$
is bounded."
($V_i$ itself is allowed to have linear growth.)
Since Bailleul \cite{be} recently solved RDEs with such coefficients,
it might be possible to extend our theorem 
to include such a case 
by just combining existing methods. 
\\
\noindent
{\bf (iv)}~
In a very recent survey \cite{bo_svy},
many results on various kinds of short time asymptotic problems
for RDEs (or Young ODE) driven by fBm are reviewed. 
For instance, 
Varadhan's estimate, which is short time asymptotics of $\log p(t, a, a')$,
was shown in \cite{boz1}
under the uniform ellipticity condition on the coefficient vector fields
when $H >1/4$.

\end{re}


%
%
\section{Moment estimate for Taylor expansion of Lyons-It\^o map}
\label{sec_two}

Let $p \in [2,3)$ be the roughness constant and let $q \in [1,2)$ be such that $1/p +1/q >1$.
We denote by $G\Omega_p ({\mathbb R}^d)$ the geometric rough path space with $p$-variation topology.
In this paper, the time interval is always $[0,1]$.
For the definition and basic properties of geometric rough paths, 
see Lyons and Qian \cite{lq}, or Lyons, Caruana, and L\'evy \cite{lcl}.

Assume that $\sigma : {\mathbb R}^n \to {\rm Mat}(n,d)$ 
and $b : [0,1] \times {\mathbb R}^n \to {\rm Mat}(n,e)$
are $C_b^{\infty}$.
For $\ve \in [0,1]$, ${\bf x} \in G\Omega_p ({\mathbb R}^d)$ and 
$h \in C_0^{q -var} ([0,1], {\mathbb R}^e)$, 
we consider the following RDE driven by the Young pairing $(\ve {\bf x}, {\bf h}) 
\in G\Omega_p ({\mathbb R}^{d+e})$;
\begin{equation}\label{rde1.eq}
d y^{\ve}_t = \sigma (y^{\ve}_t) \ve d x_t + b(\ve,  y^{\ve}_t) d h_t
\qquad
\mbox{ with }
\quad
 y^{\ve}_0 = a \in {\mathbb R}^n.
\end{equation}

It was shown in Inahama \cite{ina10} (or Inahama-Kawabi \cite{ik}) 
that the first level path of the solution 
admits a Taylor-like expansion in the deterministic sense as $\ve \searrow 0$.
Roughly speaking, the aim of this section is to prove that 
the expansion holds still true in $L^r$-sense for any $r \in [1, \infty)$, 
when ${\bf x}$ is the natural lift of fBm with $H \in (1/3, 1/2]$
or a similar Gaussian process.

We remark that the following RDE is a special case of (\ref{rde1.eq}) above:
\begin{equation}\label{rde1'.eq}
d y^{\ve}_t = \sigma ( y^{\ve}_t) ( \ve d x_t + d k_t)
+ \hat{b} (\ve, y^{\ve}_t) d \lambda_t
\qquad
\mbox{ with }
\quad
y^{\ve}_0 = a \in {\mathbb R}^n.
\end{equation}
Here, $\sigma$ and ${\bf x}$ are as above, $\hat{b} : [0,1] \times {\mathbb R}^n \to {\mathbb R}^n$,
$\lambda_t =t$, $k \in C_0^{q -var} ( [0,1], {\mathbb R}^d)$. 
We can easily check this by setting $e=d+1$, $h =(k, \lambda)$, 
and $b = [\sigma | \hat{b}]$ (an $n \times (d+1)$ block matrix).
This type 
of RDE appears when we make a Young translation of a given RDE
driven by a scaled Gaussian rough path.

\subsection{Notations}

In this paper we work in Lyons' original framework of rough path theory. 
We borrow most of notations and terminologies from \cite{lq, lcl}.
Before we start detailed discussions, however, we need to set some additional notations.

We denote by ${\bf x}= ({\bf x}^1,{\bf x}^2)$
a generic element in $G\Omega_p ({\mathbb R}^d)$
and we write $x_t := {\bf x}^1_{0,t}$ as usual.
Conversely,  for $x \in C_0^{\al-var} ([0,1], {\mathbb R}^d)$ with $\al \in [1,2)$,
we denote the natural lift of $x$ (i.e., the smooth rough path lying above $x$)
 by the corresponding boldface letter ${\bf x}$.

Note that, for $x \in C_0^{\al-var} ([0,1], {\mathbb R}^d)$ and 
$y \in C_0^{\al-var} ([0,1], {\mathbb R}^e)$, 
$({\bf x}, {\bf y}) \in G\Omega_p ({\mathbb R}^{d+e})$ stands for the natural lift of $(x,y)$, 
not for the pair $({\bf x}, {\bf y})  \in G\Omega_p ({\mathbb R}^{d}) \times G\Omega_p ({\mathbb R}^{e})$.
In a similar way, for 
${\bf x} \in G\Omega_p ({\mathbb R}^{d})$ and 
$h \in C_0^{q -var} ({\mathbb R}^e)$ with $1/p +1/q >1$,
$({\bf x}, {\bf h}) 
\in G\Omega_p ({\mathbb R}^{d+e})$ stands for the Young pairing.
These notations may be somewhat misleading.
But, they  make many operations intuitively clear and easy to understand 
when we treat rough paths over a direct sum of many vector spaces.

For a control function $\omega$ in the sense of p. 16, \cite{lq}, 
we write $\bar{\omega} := \omega (0,1)$.
For any ${\bf x} \in G\Omega_p ({\mathbb R}^{d})$, 
\begin{equation}\label{nat_control.def}
\omega_{{\bf x}} (s,t) := 
%
\|{\bf x}^1 \|_{p -var, [s,t]}^{p}
+
\|{\bf x}^2 \|_{p/2-var, [s,t]}^{p/2}
\qquad\qquad
(0 \le s \le t \le 1)
 \end{equation}
defines a control function. 
Here, the norm on the right hand side denoted the $p/j$-variation ($j=1,2$)
restricted on the subinterval $[s,t]$.
(This control function
 is equivalent to  the one defined by Carnot-Carath\'eodory metric.)
Similarly, 
we set $\omega_{\lambda} (s,t) :=  \|\ld \|^q_{q-var, [s,t]}$ 
for $\ld \in C_0^{q-var} ([0,1], {\mathbb R}^e)$.

For $\alpha >0$ and ${\bf x} \in G\Omega_p ({\mathbb R}^{d})$,
%
%
%
%
%
%
%
set $\tau_0 (\alpha) =0$ and  
$$\tau_{i+1}  (\alpha)
=\inf \{ t \in (\tau_i (\alpha)  ,1]
~|~ 
\omega_{{\bf x}} ( \tau_i (\alpha)  ,t)    
\ge \alpha\} \wedge 1
\qquad
(i=1,2,\ldots)
$$
Define 
\begin{equation}
\label{def.locvar2}
N_{\alpha} ({\bf x}) 
=
\sup \{ i \in {\mathbb N} ~|~ \tau_i (\alpha)  <1 \}.
\end{equation} 
Superadditivity of $\omega_{{\bf x}} $ yields
$\alpha N_{\alpha} ({\bf x})  \le \overline{ \omega_{{\bf x}} }$.
This quantity (\ref{def.locvar2}) was first studied by  Cass, Litterer, and Lyons \cite{cll}.
%
%

Let $x(k)$ be continuous paths which takes values in ${\mathbb R}^{d_k}$ for $k=1,\ldots, m$.
Then, we write 
\begin{equation}\label{def.calI}
{\cal I} [x(1), \ldots, x(m)]_{s,t} = \int_{s < t_1 <\cdots < t_m <t} 
dx(1)_{t_1} \otimes \cdots  \otimes  dx(m)_{t_m}
\end{equation}
whenever the iterated integral on the right hand side makes sense.
For example, if $({\bf x}, {\bf y})$ is a smooth rough path lying above $(x,y)$, 
then its second level path is given by $({\bf x}^2,  {\cal I} [x,y], {\cal I} [y,x], {\bf y}^2)$
Slightly abusing notations, we denote by ${\cal I} [x,y]$ the
$"(x,y)$-component" of the second level path of any ${\bf z} = "({\bf x}, {\bf y})" \in 
G\Omega_p ({\mathbb R}^{d} \oplus {\mathbb R}^{e})$,
when no confusion may occur.

For brevity we will often write ${\cal V} = {\mathbb R}^d$, 
$\hat{\cal V} = {\mathbb R}^e$,
and ${\cal W} = {\mathbb R}^n$ in this section.

\subsection{ODEs for ordinary Taylor terms}

Ordinary terms in the Taylor expansion are known to satisfy a very simple ODE.
In this section we recall them, following \cite{ina10}, etc. 
We will first calculate in $1$-variational setting (i.e., the Riemann-Stieltjes sense).
After that we will continuously extend these objects to the rough path setting.

The ODE that corresponds to (\ref{rde1.eq}) is the following;
\begin{equation}\label{ode1.eq}
d y^{\ve}_t = \sigma (y^{\ve}_t) \ve d x_t + b(\ve, y^{\ve}_t) dh_t
\qquad
\mbox{ with }
\quad
 y^{\ve}_0 = a.
\end{equation}
Here, $(x,h) \in C_0^{1 -var} ([0,1], {\mathbb R}^{d+e})$.
By setting $\ve =0$,
we can easily see that the $0$th term  $\phi^0 = \phi^0 (h)$
satisfies the following ODE;
\begin{equation}\label{ode_phi0.eq}
d \phi^{0}_t =  b(0, \phi^{0}_t) dh_t
\qquad
\mbox{ with }
\quad
 \phi^{0}_0 = a.  
\end{equation}
%


ODEs for $\phi^1$ and $\phi^2$ are given as follows;
\begin{equation}\label{ode_phi1.eq}
d \phi^1_t - \nabla b (0,  \phi^0_t)  \la \phi^1_t, dh_t \ra
=
\sigma (\phi^0_t) dx_t + \partial_{\ve}  b (0, \phi^0_t) dh_t
\quad
\mbox{ with }
\quad
 \phi^{1}_0 = 0,
 \end{equation}
and 
\begin{eqnarray}
d \phi^2_t - \nabla b (0,  \phi^0_t)  \la \phi^2_t, dh_t \ra
=
\nabla \sigma (\phi^0_t)\la \phi^1_t, dx_t \ra 
+
\frac12 \nabla^2 b (0, \phi^0_t) \la \phi^1_t, \phi^1_t,  dh_t\ra
\nn\\
+
 \partial_{\ve} \nabla b (0, \phi^0_t) \la  \phi^1_t,  dh_t\ra
+ 
\frac12 \partial_{\ve}^2  b (0, \phi^0_t) dh_t
\quad
\mbox{ with }
\quad
 \phi^{2}_0 = 0.
 \label{ode_phi2.eq}
  \end{eqnarray}

ODEs for $\phi^k = \phi^k (x,h) ~(k=2, 3,4,\ldots)$ are given as follows.
A heuristic explanation for how to derive these ODEs was given in \cite{ina10}.
We write $\partial_{\ve} b$ for the partial derivative in $\ve$ 
and $\nabla b$ for the (partial) gradient in $y$ for fixed $\ve$.
\begin{equation}\label{ode_phik.eq}
d \phi^k_t - \nabla b (0,  \phi^0_t)  \la \phi^k_t, dh_t \ra
=
dA^k_t 
+d B^k_t 
\quad
\mbox{ with }
\quad
 \phi^{k}_0 = 0,
\end{equation}
where 
\begin{equation}
dA^k_t [x, h, \phi^0, \ldots, \phi^{k-1}] =
\sum_{j=1}^{k-1}  
\sum_{ i_1 + \cdots + i_j = k-1}
\frac{1}{j!}
\nabla^j \sigma (\phi^0_t) \la \phi^{i_1}_t, \ldots, \phi^{i_j}_t ,dx_t  \ra
\label{def.Ak}
\end{equation}
and 
\begin{eqnarray}
dB^k_t [x,h, \phi^0, \ldots, \phi^{k-1}] &=&
\sum_{j=2}^{k}  
\sum_{ i_1 + \cdots + i_j = k}
\frac{1}{j!}
\nabla^j b (0, \phi^0_t) \la \phi^{i_1}_t, \ldots, \phi^{i_j}_t ,dh_t  \ra
\nn\\
&&
+ 
\sum_{m=1}^{k-1}  \sum_{j=1}^{k-m}  
\sum_{ i_1 + \cdots + i_j = k-m}
\frac{1}{m! j!}
 \partial_{\ve}^m \nabla^j b (0, \phi^0_t) \la \phi^{i_1}_t, \ldots, \phi^{i_j}_t ,dh_t  \ra
\nn\\
&&
+  
\frac{1}{k!}  \partial_{\ve}^k  b (0, \phi^0_t) dh_t.
\label{def.Bk}
\end{eqnarray}
Note that in the definition of $A^k$,
the summation is taken  over all positive $i_1, \ldots, i_j$ such that $i_1 + \cdots + i_j = k-1$.
A similar remark goes for the summations in the definition of $B^k$.
(As usual we set $A^k_0 =B^k_0 =0$.)

Let us recall that
we can obtain $\phi^k$ by the variation of constants formula 
since the right hand side of (\ref{ode_phi1.eq})--(\ref{ode_phik.eq}) is known.
Set $K_t =K_t[h] = \int_0^{\cdot} \nabla b (0,  \phi^0_t)  \la \,\cdot\,, dh_t \ra$
and consider the following ${\rm Mat}(n,n)$-valued ODE;
\begin{equation}\label{odeM.eq}
dM_t = (dK_t) \cdot M_t \quad
\mbox{ with }
\quad
M_0 = {\rm Id}_n.
\end{equation}
It is easy to see that $M_t^{-1}$ exists and satisfies a similar ODE.
Using this, we can easily check that $\phi^k$ has the following expression;
\begin{equation}\label{duham.eq}
\phi^k_t = M_t \int_0^{t} M_s^{-1} dZ^k_s 
= 
Z^k_t -  M_t \int_0^{t} dM_s^{-1} \cdot Z^k_s.
\end{equation}
Here, $Z^k_t$ (with $Z^k_0=0$) is a shorthand for
the right hand side of (\ref{ode_phi1.eq})--(\ref{ode_phik.eq}). 
Finally, we set 
\begin{equation}\label{def.r}
r_{\ve}^{k+1} = y^{\ve} - (\phi^0  + \ve \phi^1 +\cdots + \ve^k \phi^k ).
\end{equation}
It is obvious that for each $\ve \in [0,1]$ and $k \in {\mathbb N}$
\begin{equation}\label{def.map.1var}
(x,h) \mapsto (x,h, y^{\ve}, \phi^0, \ldots, \phi^k, r_{\ve}^{k+1} )
\end{equation}
is continuous from $C_0^{1-var}([0,1], {\cal V} \oplus \hat{\cal V})$
to 
$C^{1-var}([0,1], {\cal V} \oplus \hat{\cal V} \oplus {\cal W}^{ \oplus k+3})$.
%
%
%
%
%
It is known that this map extends to a continuous map 
with respect to the rough path topology in the following sense
(after the initial values are suitably adjusted, precisely speaking. 
Note that $y^{\ve}_0 =a= \phi^0_0$.) 
\begin{pr}\label{pr.map.rp}
Let $2 \le p <3$ and $1 \le q <2$ such that $1/p +1/q >1$.
Then, for each $\ve \in [0,1]$ and $k \in {\mathbb N}$,
the map  (\ref{def.map.1var}) naturally extends to 
the following locally Lipschitz continuous map;
\begin{equation}
G\Omega_p ({\cal V} ) \times C_0^{q-var}(\hat{\cal V} ) 
\ni
({\bf x},h) \mapsto ({\bf x}, {\bf h}, {\bf y}^{\ve}, \bm{\phi}^0, \ldots, \bm{\phi}^k, {\bf r}_{\ve}^{k+1} )
\in
G\Omega_p ({\cal V} \oplus \hat{\cal V} \oplus {\cal W}^{ \oplus k+3}).
\nn
\end{equation}
\end{pr}

\Proof
This was already shown in \cite{ina10} for arbitrary $p \ge 2$.
Here, we only give a sketch of proof for later use.

First, $({\bf x}, {\bf h})$ is just Young pairing of ${\bf x}$ and $h$.
Since $({\bf y}^{\ve}, \bm{\phi}^0)$ 
is a unique solution of an RDE driven by $({\bf x}, {\bf h})$,
we obtain 
$({\bf x}, {\bf h}, {\bf y}^{\ve}, \bm{\phi}^0)$.
Next, assume that we have 
$({\bf x}, {\bf h}, {\bf y}^{\ve}, \bm{\phi}^0, \ldots, \bm{\phi}^{k-1})$.
Then, 
$A_k +B_k$ on the right hand side of (\ref{ode_phik.eq})
can be interpreted 
as a rough path integral, we obtain 
$({\bf x}, {\bf h}, {\bf y}^{\ve}, \bm{\phi}^0, \ldots, \bm{\phi}^{k-1}, A_k +B_k)$.
For $\tilde{M}_t := {\rm Id}_{{\cal V} \oplus \hat{\cal V} \oplus {\cal W}^{ \oplus k+1}} \oplus M_t$,
we can use a rough path version of variation of constant method
to obtain 
$({\bf x}, {\bf h}, {\bf y}^{\ve}, \bm{\phi}^0, \ldots, \bm{\phi}^{k})$.
(Observe (\ref{duham.eq}) above.)
Finally, since $r^{k+1}_{\ve}$ is a linear combination of 
$y^{\ve}, \phi^0, \ldots, \phi^{k}$,
we obtain 
$({\bf x}, {\bf h}, {\bf y}^{\ve}, \bm{\phi}^0, 
\ldots, \bm{\phi}^k, {\bf r}_{\ve}^{k+1} )$.
\QED

By the following proposition, this expansion can be called 
a Taylor(-like) expansion of  Lyons-It\^o map.

\begin{pr}\label{pr.map.rp2}
Keep the same notations and assumptions as in Proposition \ref{pr.map.rp} above.
Then, the following {\bf (i)} and {\bf (ii)} hold.
\\
\noindent
{\bf (i)}~
For any any $\rho >0$ and $k=1,2,\ldots$, there exists a positive constants $C=C(\rho, k)$ 
which satisfies that
\[
\| (\bm{\phi}^k)^1 \|_{p -var} \le C (1+ \overline{\omega_{ {\bf x}}}^{1/p})^k.
\]
for any ${\bf x} \in G\Omega_p ({\cal V} )$ and 
any $h \in C_0^{q-var}([0,1], \hat{\cal V} )$ with $\|h\|_{q -var} \le \rho$..
\\
\noindent
{\bf (ii)}~
For any $\rho_1, \rho_2 >0$ and $k=1,2,\ldots$,
there exists a positive constants $\tilde{C}= \tilde{C}(\rho_1, \rho_2, k)$, which is independent of $\ve$
and satisfies that
\[
\| ({\bf r}^{k+1}_{\ve})^1 \|_{p -var} \le \tilde{C} (\ve + \ve \overline{\omega_{  {\bf x}}}^{1/p}) ^{k+1} 
\]
for any ${\bf x} \in G\Omega_p ({\cal V} ) $ with 
$\overline{\omega_{ \ve {\bf x}}}^{1/p}   =  \ve \overline{\omega_{ {\bf x}}}^{1/p}   \le \rho_1$
and any $h \in C_0^{q-var}([0,1], \hat{\cal V} )$ with $\|h\|_{q -var} \le \rho_2$.
\end{pr}

\Proof 
This was already shown in \cite{ina10} for arbitrary $p \ge 2$.
In that paper,  estimates not only for the first level path,
but also for the higher level paths are given.
\QED

\begin{re}
In a very recent preprint \cite{be2}, 
Bailleul gave a simplified proof of Propositions \ref{pr.map.rp} and \ref{pr.map.rp2}
for any $p \ge 2$
in the framework of Gubinelli's controlled path theory.
\end{re}

\subsection{Main results in this section}
In this subsection 
we state the main result of this section, 
that is, moment estimates for Taylor expansion of Lyons-It\^o map.
%
%
We will prove this theorem rigorously in subsequent subsections.
Note that $\eta_k$ may depend on $k,{\bf x}, h, p, q$,
but not on $\ve$.

%

\begin{tm}\label{thm.mom_r}
Let $2 \le p < 3$ and $1 \le q <2$ such that $1/p + 1/q >1$
and
let $h \in C_0^{q -var} ([0,1], {\mathbb R}^e)$.
Assume that 
${\bf x}$ is a $G\Omega_p ({\mathbb R}^d )$-valued random variable 
such that 
{\rm (a)}~ $\overline{\omega_{{\bf x}}}= \omega_{{\bf x}} (0,1) \in \cap_{1 \le r <\infty} L^r$ 
and 
{\rm (b)}~ 
$\exp ( N_{\alpha} ({\bf x}) )
 \in  \cap_{1 \le r <\infty} L^r$ for any $\alpha >0$ .

Then, for any ${\bf x}$, $h$,  and $k \in {\mathbb N}$,
there exist control functions $\eta_k = \eta_{k,{\bf x}, h }$
such that the following {\rm (i)}-- {\rm (iii)} hold:
\\
\noindent
{\rm (i)}~$\eta_{k}$ are non-decreasing in $k$, i.e., 
$\eta_{k,{\bf x}, h}(s,t) \le \eta_{k+1,{\bf x}, h} (s,t)$ 
for all $k, {\bf x}, h, (s,t)$.
\\
\noindent
{\rm (ii)}~$\overline{\eta_{k,{\bf x}, h} } \in  \cap_{1 \le r <\infty} L^r$ 
for all $k, h$.
\\
\noindent
{\rm (iii)}~
For all $\ve \in (0,1]$, $k\in {\mathbb N}$, $h$, ${\bf x}$, and $0 \le s \le t \le t$, $j=1, 2$,  we have 
\[
\bigl|  
\bigl(
{\bf x}, {\bf h}, {\bf y}^{\ve}, \bm{\phi}^0, \ldots, \bm{\phi}^k, \ve^{-(k+1)} {\bf r}_{\ve}^{k+1}
\bigr)^j_{s,t}
\bigr| \le \eta_{k,{\bf x}, h} (s,t)^{j/p}.
\]

In particular, for all $k\in {\mathbb N}$ and $h$, 
$\| (\bm{\phi}^k)^1\|_{p -var} \in \cap_{1 \le r <\infty} L^r$
and 
$ \| ({\bf r}_{\ve}^{k+1})^1 \|_{p -var} = O (\ve^{k+1})$ in $L^r$ for any $1 <r <\infty$.
\end{tm}

\begin{re}\label{re.int.fbrp}
{\rm (1)}~
Examples of Gaussian processes whose 
rough path lifts satisfy the integrability assumptions 
\[
\omega_{{\bf x}} (0,1) \in \cap_{1 \le r <\infty} L^r
\quad 
\mbox{ and }
\quad
\exp ( N_{\alpha} ({\bf x})  )
 \in  \cap_{1 \le r <\infty} L^r
 \quad (\forall \alpha >0)
 \]
can be found in Friz and Oberhauser \cite{fo} (a Fernique-type theorem)
and 
Cass, Litterer, and Lyons \cite{cll} (Integrability of $N_{\alpha}$).
FBm with Hurst parameter $H \in (1/4, 1/2]$ is a typical example. 
%
\\
{\rm (2)}~
The estimate above is actually uniform in $h$
when it varies in a bounded subset in $q$-variation space.
(But, the uniform version is not needed in this paper.)
\end{re}

\subsection{Proof of Theorem \ref{thm.mom_r} for $k=0$}

The rest of this section is devoted to showing Theorem \ref{thm.mom_r}.
Without loss of generality we may assume that the initial value $a=0$.
In this proof
$c_1, c_2, \ldots$ stands for unimportant positive constants,
which is independent of $\ve \in (0,1]$ and ${\bf x}$, 
but may depend on $p, q, \sigma, b, \|h\|_{q -var}$, etc.
%
We say that a geometric $p$-rough path ${\bf x}$ is controlled by 
a control function $\omega$ 
if $|{\bf x}^j_{s,t} | \le \omega(s,t)^{j/p}$
for any $s \le t$ and $j=1,2$.

The expansion of the It\^o map 
in the deterministic case is already given in \cite{ina10, ik} by mathematical induction. 
We will closely look at it and check the integrability holds or not.
In this subsection,
we will obtain the moment estimates of ${\bf r}_{\ve}^{1}$.
Surprisingly, 
for those who understand the proof for the deterministic sense, 
the most difficult part is this initial step of the induction.
However, that problem is somewhat similar  to the moment estimates of 
Jacobian process driven by Gaussian rough paths,
which was solved by Cass, Litterer, and Lyons \cite{cll}.
In the sequel we will check that their method also applies to this kind of problem 
as they conjectured in \cite{cll}.

Now we prove Theorem \ref{thm.mom_r} for $k=0$.
Set $\omega_{h} (s,t)= \|h\|_{q -var, [s,t]}^q$
and 
$\omega_{{\bf x},h} (s,t) = \omega_{{\bf x}} (s,t) + \omega_{h} (s,t)$.
Then, the Young pairing $({\bf x}, {\bf h}) \in G\Omega_p ( {\cal V} \oplus \hat{\cal V})$
is controlled by $c_1 (1+ \overline{\omega_{{\bf x},h} })^{c_2} \omega_{{\bf x},h} $, that is, 
\[
\bigl| ({\bf x}, {\bf h})^i_{s,t} \bigr| \le  \bigl\{ c_1 (1+ \overline{\omega_{{\bf x},h} })^{c_2} \omega_{{\bf x},h} (s,t) \bigr\}^{i/p}
\]
for all $i$ and $(s,t)$.

Next we consider 
$( {\bf y}^{\ve}, \bm{\phi}^0)$ which is a solution of a ${\cal W}^{\oplus 2}$-valued RDE driven by 
$({\bf x}, {\bf h})$.
Since the $C_b^{[p]+1}$-norm of the coefficients of the RDE is bounded in $\ve$,
$({\bf x}, {\bf h},  {\bf y}^{\ve}, \bm{\phi}^0) 
\in G\Omega_p ( {\cal V} \oplus \hat{\cal V} \oplus {\cal W}^{\oplus 2})$
is controlled by 
$c_3 (1+ \overline{\omega_{{\bf x},h} })^{c_4} \omega_{{\bf x},h}$.

It is easy to see
from (\ref{ode1.eq}) and (\ref{ode_phi0.eq}) that $r^{1}_{\ve, t}$ satisfies 
the following equation in the $1$-variational setting;
\begin{equation}\label{ode_r1.eq}
\frac{1}{\ve} d r^{1}_{\ve, t}
=
\sigma (y^{\ve}_t)  d x_t + \frac{1}{\ve} \{ b(\ve, y^{\ve}_t)  -  b(0, \phi^{0}_t)   \}dh_t
\qquad
\mbox{ with }
\quad
r^{1}_{\ve, 0}=0.
\end{equation}
The first term on the right hand side can be interpreted as 
a rough path integration of a $C^{[p]+1}_b$ one-form
along $({\bf x}, {\bf h},  {\bf y}^{\ve}, \bm{\phi}^0)$.
Hence, 
$({\bf x}, {\bf h},  {\bf y}^{\ve}, \bm{\phi}^0, \int \sigma ( {\bf y}^{\ve}) d{\bf x} ) $
is controlled by 
$c_5 (1+ \overline{\omega_{{\bf x},h} })^{c_6} \omega_{{\bf x},h}$, namely, 
\begin{equation}\label{r1l.ineq}
\bigl| \bigl( {\bf x}, {\bf h},   {\bf y}^{\ve}, 
\bm{\phi}^0, \int \sigma ( {\bf y}^{\ve}) d{\bf x}  \bigr)^i_{s,t}  \bigr| 
\le  \bigl\{ c_5 (1+ \overline{\omega_{{\bf x},h} })^{c_6} \omega_{{\bf x},h} (s,t) \bigr\}^{i/p}
\end{equation}
for all $i$ and $(s,t)$.
With (\ref{r1l.ineq}) in hand we have only to obtain 
a nice estimate of $q$-variation norm of the second term  on the right hand side of (\ref{ode_r1.eq}).

Let us  estimate the first level path of ${\bf r}_{\ve}^{1}$,
that is the difference of the first level paths of
${\bf y}^{\ve}$ and $\bm{\phi}^0$.
${\bf y}^{\ve}$ 
and $\bm{\phi}^0$ are the solutions of 
the RDEs (whose coefficient are the $n \times (d+e)$-matrix
$[\sigma, b(\ve, \,\cdot\,)]$ and $[\sigma, b(0, \,\cdot\,)]$, resp.)
driven by $(\ve {\bf x}, {\bf h})$ and $({\bf 0}, {\bf h})$, resp.
A useful estimate of 
difference of two solutions of RDEs can be found in 
Theorem 10.26, pp. 233--236, \cite{fvbk}.
({\bf 0} denotes a rough path such that ${\bf x}^j \equiv 0$ for $j=1,2$.)
There are positive constants $c_7, c_8$ such that 
$C_b^{[p]}$-norm of $[\sigma, b(\ve, \,\cdot\,)]$ and 
$[\sigma, b(\ve, \,\cdot\,)]- [\sigma, b(0, \,\cdot\,)]$
are dominated by $c_7$ and $c_8 \ve$, respectively.
Set $\omega' = c_9 \omega_{{\bf x},h} (s,t) $.
If we take $c_9$ sufficiently large, we have 
$\|(\ve {\bf x}, {\bf h}) \|_{p -\omega'} \le 1$ for all $\ve \in [0,1]$
(see Chapter 8, \cite{fvbk} the definition and $\| \,\cdot\, \|_{p -\omega'}$ 
and details)
and
$(\ve {\bf x}, {\bf h})$ satisfies the assumption {\rm (iii)}, Theorem 10.26, \cite{fvbk}.
Note also that 
$
\bigl| (\ve {\bf x}, {\bf h})^j_{s,t}  - ( {\bf 0}, {\bf h})^j_{s,t} \bigr|
\le 
\ve 
c_{10} (1+ \omega' (S,T) )^{c_{11}} \omega' (s,t)^{j/p}$ for any 
$S, T, s, t$
with $S \le s \le t \le T$.

Then, (a trivial modification of) Theorem 10.26, \cite{fvbk} implies that,
on any subinterval $[S,T] \subset [0,1]$, there exists a constant $c_{12}>0$ such that
\begin{eqnarray}
\lefteqn{
\bigl| ({\bf y}^{\ve})^1_{s,t} - (\bm{\phi}^0)^1_{s,t} \bigr|
}
\nn\\
&\le&
c_{12} \bigl[
c_7 | ({\bf y}^{\ve})^1_{0,S} - (\bm{\phi}^0)^1_{0,S}| + c_8 \ve 
+
 \ve 
c_{10} (1+  \omega' (S,T) )^{c_{11}} \bigr]  \omega' (s,t)^{1/p}
\exp (c_{12} c_7^p \omega' (S,T))
\nn\\
&\le&
c_{13}
\bigl[
| ({\bf y}^{\ve})^1_{0,S} - (\bm{\phi}^0)^1_{0,S}| + 
\ve 
 (1+ \omega_{{\bf x},h }(S,T) )^{c_{13}}
\bigr] 
 \omega_{{\bf x},h} (s,t)^{1/p}
\exp (c_{13} \omega_{{\bf x},h}(S,T) )
\label{FVthm1026}
\end{eqnarray}
for any $S \le s \le t \le T$.

Let $\tau_{i} =\tau_{i} (\alpha)$ be as in the definition of $N_{\alpha} ({\bf x}) $ in (\ref{def.locvar2}).
We choose $\alpha>0$ so small that
$c_{13} (1+ 2\alpha)^{ c_{13}}  (2\alpha)^{1/p} \le 1$ holds.
Consider each subinterval $I_i :=[\tau_{i-1}, \tau_{i}]$ ($i=1,2,\ldots, N_{\alpha} ({\bf x}) $).
Let 
$\{\tau_{i-1}= \sigma^{(i)}_0 < \sigma^{(i)}_1 <\cdots< \sigma^{(i)}_{K_i}  =\tau_{i} \}$
be a partition of $I_i$ 
such that 
$\omega_h (\sigma^{(i)}_{j-1}, \sigma^{(i)}_j) =\alpha$ for $1 \le j \le K_i-1$
and 
$\omega_h (\sigma^{(i)}_{K_i-1}, \sigma^{(i)}_{K_i}) \le \alpha$.
It is easy to see that $K_i -1 \le \omega_h ( \tau_{i-1}, \tau_{i}) /\alpha$.
Let $\{0=t_0 <t_1 <\cdots < t_J =1\}$ be all $\sigma^{(i)}_{j}$'s in increasing order.
The total number $J$ of 
the subintervals is now at most
\[
J= \sum_{i=1}^{ N_{\alpha} ({\bf x}) +1}  K_i 
\le
 \sum_{i=1}^{ N_{\alpha} ({\bf x}) +1} 
 \bigl( 
  1+  \frac{\omega_h ( \tau_{i-1}, \tau_{i}) }{\alpha}
     \bigr)
     \le 
       N_{\alpha} ({\bf x}) +1 +  \frac{\|h\|^q_{q -var} }{\alpha}.
                \]

On each subinterval $\hat{I}_i :=[t_{i-1}, t_i]$, $\omega_{{\bf x},h} ( t_{i-1}, t_i) \le 2\alpha$.
Hence we have from (\ref{FVthm1026}) that 
\[
\bigl| ({\bf y}^{\ve})^1_{s,t} - (\bm{\phi}^0)^1_{s,t} \bigr|
\le 
\bigl[
| ({\bf y}^{\ve})^1_{0, t_{i-1}} - (\bm{\phi}^0)^1_{0, t_{i-1}}| + 
\ve 
\bigr] 
\exp (c_{13} \omega_{{\bf x},h}( t_{i-1}, t_i ))
\]
for any $t_{i-1} \le s \le t \le t_{i}$.
By mathematical induction, we have 
\begin{eqnarray*}
| ({\bf y}^{\ve})^1_{0, t_{i-1}} - (\bm{\phi}^0)^1_{0, t_{i-1}}|
&\le&
\ve
\prod_{k=1}^{i-1}  \bigl\{
1 + \exp (c_{13} \omega_{{\bf x},h}( t_{k-1}, t_k ))
\bigr\} 
\nn\\
&\le& 
\ve
2^{i-1}  \exp \bigl(c_{13} \sum_{k=1}^{i-1} \omega_{{\bf x},h}( t_{k-1}, t_k )  \bigr).
\end{eqnarray*}
Putting this back into (\ref{FVthm1026}), we have on each interval $\hat{I}_i$, 
\begin{eqnarray*}
\bigl| ({\bf y}^{\ve})^1_{s,t} - (\bm{\phi}^0)^1_{s,t} \bigr|
&\le&
\ve 
\Bigl\{
c_{13} 2^{i-1} \exp \bigl(c_{13} \sum_{k=1}^{i-1} \omega_{{\bf x},h}( t_{k-1}, t_k )  \bigr)
+
(2\alpha)^{-1/p}
\Bigr\} 
 \nn\\
 && 
 \times 
 \omega_{{\bf x},h} (s,t)^{1/p}
  \exp (c_{13} \omega_{{\bf x},h}( t_{i-1}, t_i ))
   \nn\\
    &\le& 
  \ve      c_{14} 
  \exp \bigl( J \log 2 + c_{13} \sum_{k=1}^{J} \omega_{{\bf x},h}( t_{k-1}, t_k )  \bigr)
   \omega_{{\bf x},h} (s,t)^{1/p}
             \nn\\
             &\le&
               \ve      c_{14} 
  \exp \Bigl[
      ( N_{\alpha} ({\bf x}) +1 +  \frac{\|h\|^q_{q -var} }{\alpha}) \log 2 
      \nn\\
      &&
     \qquad \qquad  \qquad + 
       c_{13} \{ \alpha (N_{\alpha} ({\bf x}) +1)+  \|h\|^q_{q -var} \}
                   \Bigr]
   \omega_{{\bf x},h} (s,t)^{1/p}
      \nn\\
      &\le& 
          \ve      c_{15}  
          \exp (c_{16} N_{\alpha} ({\bf x}) ) 
            \omega_{{\bf x},h} (s,t)^{1/p}.
                               \end{eqnarray*}
Here, the positive constants $c_i~(14 \le i \le 16)$ depend on $\alpha$, too.
Since there are $J$ subintervals, we have on the whole interval that
\begin{eqnarray}
\bigl| ({\bf y}^{\ve})^1_{s,t} - (\bm{\phi}^0)^1_{s,t} \bigr|
&\le&
J^{1 -1/p} \ve      c_{15}  
          \exp (c_{16} N_{\alpha} ({\bf x}) ) 
            \omega_{{\bf x},h} (s,t)^{1/p}
            \nn\\
              &\le&
              \ve      c_{17}
              \bigl(
               N_{\alpha} ({\bf x}) +1  
                               \bigr)^{1 -1/p}  
                               \exp (c_{16} N_{\alpha} ({\bf x}) ) 
                                \omega_{{\bf x},h} (s,t)^{1/p}
                                \nn\\
                                &\le& 
                                  \ve c_{18}\exp (c_{18} N_{\alpha} ({\bf x}) ) 
                                \omega_{{\bf x},h} (s,t)^{1/p}
                  \label{pvy-phi.eq}
                  \end{eqnarray}
for any $0 \le s \le t \le 1$.  This is the most difficult part in this subsection.
For brevity we set a control function $\xi_1$ by 
$\xi_1 (s,t)^{1/p} = c_{18}\exp (c_{18} N_{\alpha} ({\bf x}) ) 
                                \omega_{{\bf x},h} (s,t)^{1/p}$. 
%
%
Obviously, 
$\overline{ \xi_1} \in \cap_{1 <r<\infty} L^r$ by assumption.

%
We see from (\ref{r1l.ineq}) and (\ref{pvy-phi.eq}) that
\begin{eqnarray*}
\lefteqn{
\Bigl|
 \{ b(\ve, y^{\ve}_t)  -  b(0, \phi^{0}_t)   \}
 -
  \{ b(\ve, y^{\ve}_s)  -  b(0, \phi^{0}_s)   \}
   \Bigr|
}
\nn\\
&\le&
\| \nabla b\|_{\infty} |({\bf y}^{\ve})^1_{s,t}  - (\bm{\phi}^0)^1_{s,t}| 
+
\Bigl( 
\ve  \|\partial_{\ve} \nabla b\|_{\infty}
+
2\| \nabla^2 b\|_{\infty}  \|y^{\ve} - \phi^0  \|_{\infty}
\Bigr)
| (\bm{\phi}^0)^1_{s,t}| 
\nn\\
&\le&
\ve
\Bigl[
 \| \nabla b\|_{\infty} \xi_{1} (s,t)^{1/p}
+
\Bigl( 
 \|\partial_{\ve} \nabla b\|_{\infty}
+
2\| \nabla^2 b\|_{\infty}  \overline{ \xi_1}
\Bigr)
 \bigl\{ c_5 (1+ \overline{\omega_{{\bf x},h} })^{c_6} \omega_{{\bf x},h} (s,t) \bigr\}^{1/p}
\Bigr]
\nn\\
&\le&
\ve
2^{(p-1)/p}
\Bigl[
 \| \nabla b\|_{\infty}^p \xi_{1} (s,t)
+
\Bigl( 
 \|\partial_{\ve} \nabla b\|_{\infty}
+
2\| \nabla^2 b\|_{\infty}  \overline{ \xi_1}
\Bigr)^p
  c_5 (1+ \overline{\omega_{{\bf x},h} })^{c_6} \omega_{{\bf x},h} (s,t) 
\Bigr]^{1/p}.
\end{eqnarray*}
We denote the right hand side by $\ve \xi_2 (s,t)^{1/p}$.
Then, $\xi_2$ is a control function such that
$\overline{ \xi_2} \in \cap_{1 <r<\infty} L^r$.

From a basic property of Young integration and
the above estimate, we have
\begin{eqnarray}
\frac{1}{\ve} \Bigl| \int_s^t \{ b(\ve, y^{\ve}_t)  -  b(0, \phi^{0}_t)   \}dh_t \Bigr|
\le 
c_{19} 
\bigl(
1 + \overline{ \xi_2} + \overline{ \omega_h}
\bigr)^{c_{19}}
\{ \xi_2 (s,t) + \omega_h (s,t) \}^{1/q}.
\label{buri.ineq}
\end{eqnarray}
In particular, the Young integral on the left hand side above
is of finite $q$-variation.

From (\ref{ode_r1.eq}), (\ref{r1l.ineq}), (\ref{buri.ineq}) and 
a basic property of Young pairing, we have 
\begin{equation}
\nn
\bigl| \bigl( {\bf x}, {\bf h},   {\bf y}^{\ve}, 
\bm{\phi}^0, \ve^{-1} {\bf r}_{\ve}^1  \bigr)^i_{s,t}  \bigr| 
\le  \xi_3 (s,t)^{i/p}
\end{equation}
for some control function $\xi_3$ such that $\overline{ \xi_3} \in \cap_{1 <r<\infty} L^r$.
This $\xi_3$ can be written as a simple combination of 
control functions which appear 
on the right hand sides on (\ref{r1l.ineq}) and (\ref{buri.ineq})
and  is independent of $\ve$.
Thus, we have shown
Theorem \ref{thm.mom_r} for $k=0$

\subsection{Proof of Theorem \ref{thm.mom_r} for general $k \ge 1$}

In this subsection we prove  
Theorem \ref{thm.mom_r} for $k$, assuming that it holds for the cases up to $k-1$.
%
%
In the proof of the deterministic case in \cite{ina10, ik}, 
it is explained how to obtain an estimate of ${\bf r}_{\ve}^{k+1}$,
which can be expressed as a rough path integral along 
$(
{\bf x}, {\bf h}, {\bf y}^{\ve}, \bm{\phi}^0, \ldots, \bm{\phi}^{k-1})$.
%
%
%
%

Our strategy is quite simple. 
We carefully look at the proof in \cite{ina10, ik} once again and make sure that
every operation is "of at most polynomial order."  
Therefore, for those who already know the proof for the deterministic case, 
this subsection is not very difficult.
Since the full proof is quite lengthy, we only give a sketch of proof here.

Let us calculate $r^{k+1}_{\ve}$. 
From (\ref{ode1.eq})--(\ref{ode_phik.eq}), we have
\begin{eqnarray}
\lefteqn{
d r^{k+1}_{\ve,t}  - \nabla b (0, \phi^0_t)  \la  r^{k+1}_{\ve,t}, dh_t \ra
=
\Bigl[ 
\sigma (y^{\ve}_t) \ve d x_t  -\sum_{l=1}^k \ve^l dA^l_t
\Bigr]
}
\nn\\
&&
+ \Bigl[
b(\ve, y^{\ve}_t) dh_t
-b(0, \phi^{0}_t) dh_t
- \nabla b (0, \phi^0_t)  \la y^{\ve}_t -  \phi^0_t , dh_t \ra
 - \sum_{l=1}^k \ve^l dB^l_t
\Bigr]
\nn\\
&=:&
dI_t^{k+1}  + d J_t^{k+1} , \qquad  \mbox{with $r^{k+1}_{\ve,0} =0$.}
\label{ode.rk+1.eq}
\end{eqnarray}
Here, $I^{k+1} $ and $J^{k+1} $ stand for sums of the integrals 
with respect to $x$ and $h$, respectively.
Observe  the right hand side of (\ref{ode.rk+1.eq}). 
There are only $x, h, y^{\ve}, \phi^0, \ldots, \phi^{k-1}$
(and no $\phi^{k}$).
See (\ref{def.Ak}) and (\ref{def.Bk}).
Therefore, the right hand side 
can be regarded as a rough path integral along 
$({\bf x}, {\bf h}, {\bf y}^{\ve}, \bm{\phi}^0, \ldots, \bm{\phi}^{k-1}, \ve^{-k} {\bf r}^{k}_{\ve})$.
As a result we obtain
$$
({\bf x}, {\bf h}, {\bf y}^{\ve}, \bm{\phi}^0, \ldots, \bm{\phi}^{k-1}, 
 \ve^{-k} {\bf r}^{k}_{\ve}, {\bf I}^{k+1}  + {\bf J}^{k+1} )
\in G\Omega_p ({\cal V} \oplus  \hat{\cal V} \oplus {\cal W}^{ \oplus k+3}).
$$
We will prove that the rough path above is controlled by a nice control function
with moments of all order.
Note that $J^{k+1}$ is a path of finite $q$-variation 
and hence the above rough path is a Young translation of 
$({\bf x}, {\bf h}, \ldots, \ve^{-k} {\bf r}^{k}_{\ve}, {\bf I}^{k+1})$ by $J^{k+1}$.

From Taylor expansion and 
the way the rough path integral is defined, 
we can see that the above rough path satisfies 
essentially the same estimate as in Theorem \ref{thm.mom_r}, {\rm (iii)} as follows.

\begin{lm}
\label{lm.kinouII}
Keep the same notations and assumptions as in Theorem \ref{thm.mom_r}.
Assume that Theorem \ref{thm.mom_r} holds for the cases $1,2,\ldots, k-1$.
Then, there exists a control function $\xi = \xi_{{\bf x}, h} $ such that
$\eta_{k-1,{\bf x}, h}(s,t) \le \xi_{{\bf x}, h} (s,t)$, 
$\overline{\xi} \in  \cap_{1 \le r <\infty} L^r$, 
and
\begin{equation}\label{yotsu.ineq}
\bigl|
\bigl(
{\bf x}, {\bf h}, {\bf y}^{\ve}, \bm{\phi}^0, \ldots, \bm{\phi}^{k-1}, \ve^{-k} {\bf r}_{\ve}^{k},
\ve^{-(k+1)} ( {\bf I}^{k+1}+  {\bf J}^{k+1})
\bigr)^j_{s,t}
\bigr| \le \xi_{{\bf x}, h} (s,t)^{j/p}.
\end{equation}
for all $0 \le s \le t \le 1$ and $j=1,2$.
(Note that $\xi$ may not depend on $\ve$.)
\end{lm}

Let $M$ be as in (\ref{odeM.eq}).
Then, $M$ and $M^{-1}$ are deterministic, depends only on $h$,
and are of finite $q$-variation.
We see from (\ref{ode.rk+1.eq}) that at least formally  
$$r^{k+1}_{\ve,t} 
= 
M_t \int_0^{t} M_s^{-1} d [ I_s^{k+1} + J_s^{k+1}] 
=
[ I_t^{k+1} + J_t^{k+1}]
-
 M_t \int_0^{t} dM_s^{-1} \cdot [ I_s^{k+1} + J_s^{k+1}].
$$
Note that the last expression takes the form of Young translation.

To be more precise, set 
$\tilde{M}_t := {\rm Id}_{{\cal V} \oplus \hat{\cal V} \oplus {\cal W}^{ \oplus k+2}} \oplus M_t$
and
apply (a rough path version of) variation of constant method as in (\ref{duham.eq})
to the rough path in (\ref{yotsu.ineq})  in Lemma \ref{lm.kinouII} above.
Then, we obtain
$\bigl(
{\bf x}, {\bf h}, {\bf y}^{\ve}, \bm{\phi}^0, \ldots, \bm{\phi}^{k-1}, \ve^{-k} {\bf r}_{\ve}^{k},
\ve^{-(k+1)} {\bf r}_{\ve}^{k+1}
\bigr)$.
We can easily see that this rough path satisfies 
the same inequality as in (\ref{yotsu.ineq}) (if $\xi$ is suitably replaced).

Note that $\phi^{k} = \ve^{-k} r_{\ve}^{k} - \ve \{ \ve^{-(k+1)} r_{\ve}^{k+1}\}$.
By applying a simple linear map to the above rough path,  
we can obtain 
$\bigl(
{\bf x}, {\bf h}, {\bf y}^{\ve}, \bm{\phi}^0, \ldots, \bm{\phi}^{k-1}, \bm{\phi}^{k},
\ve^{-(k+1)} {\bf r}_{\ve}^{k+1}
\bigr)$.
Since the operator norm of this $\ve$-dependent linear map  is bounded in $\ve$,
this rough path also
satisfies 
the same inequality as in (\ref{yotsu.ineq}) 
(if $\xi$ is suitably replaced by another control function, 
which we call $\eta_{k,{\bf x}, h}$).
This is the sketch of proof of Theorem \ref{thm.mom_r}.
\rule[0em]{0.5ex}{1.8ex}

\subsection{Remark for fractional order case}
In this subsection we consider the case where the coefficients of RDEs are 
of fractional order in $\ve$
and present analogous results to Proposition \ref{pr.map.rp}, Proposition \ref{pr.map.rp2}, 
and Theorem \ref{thm.mom_r}.
The contents of this subsection will be used in later sections.

In this subsection we assume that $1/3<1/p <H \le 1/2$.
Let $\sigma : {\mathbb R}^n \to {\rm Mat}(n,d)$ 
and $b : {\mathbb R}^n \to {\mathbb R}^n$ be $C_b^{\infty}$.
Let ${\bf x} \in G\Omega_p  ({\mathbb R}^d)$
and $h \in C_0^{q -var} ({\mathbb R}^d)$ with $1/p +1/q >1$ and we set $\lambda_t =t$.
We consider the following RDE driven by the Young pairing $(\ve {\bf x}, {\bf h}, {\bm \lambda})$;
\begin{eqnarray}\label{rde_fr.eq}
d \tilde{ y}^{\ve}_t &=& \sigma (\tilde{ y}^{\ve}_t) ( \ve d x_t + dh_t)
+ \ve^{1/H} b (\tilde{y}^{\ve}_t) dt
\nn
\\
&=&
\sigma (\tilde{y}^{\ve}_t)  \ve dx_t + 
\bigl[
\sigma (\tilde{y}^{\ve}_t)  dh_t
+ \ve^{1/H} b (\tilde{ y}^{\ve}_t) d\lambda_t \bigr]
\qquad
\mbox{ with }
\quad
\tilde{y}^{\ve}_0 = a \in {\mathbb R}^n.
\end{eqnarray}
This is a variant of RDE (\ref{rde1'.eq}).
Strictly speaking, unless $H =1/2$ 
the results in previous subsections cannot be used for RDE (\ref{rde_fr.eq}).
With minor modifications, however, 
similar results hold in this case, too.
We will explain it below. (Proofs are essentially the same and will be omitted).

Let us fix some notations for fractional order expansions.
For 
$$ \Lambda_1 =\{ n_1 + \frac{n_2 }{H}   ~|~ n_1, n_2 \in {\mathbb N} \},$$
let 
$0=\kappa_0 < \kappa_1 <\kappa_2 <\cdots$ 
be all elements of $\Lambda_1$ in increasing order.
More concretely, leading  terms are as follows if $H \in (1/3,1/2)$;
\begin{eqnarray}
(\kappa_0, \kappa_1, \kappa_2, \ldots)
&=&
(0, \, 1,  \,  2, \, \frac{1}{H},
\, 3,  \, 1+  \frac{1}{H} , \,  4, \, 2+ \frac{1}{H} , \, 5\wedge  \frac{2}{H}, \ldots).
\label{eq.index}
\end{eqnarray}
If $H=1/2$, then $\Lambda_1 = {\mathbb N}$.

Instead of (\ref{def.r})  
the Taylor expansion of Lyons-It\^o map takes the following form; 
\begin{equation}\label{def.r.2}
r_{\ve}^{ \kappa_{k+1}} = \tilde{y}^{\ve} - (\phi^0  + \ve^{\kappa_1} \phi^{\kappa_1} +\cdots 
+ \ve^{\kappa_k} \phi^{\kappa_k} ).
\end{equation}
In this case,  $\phi^{\kappa_k}$ is the term of "order $\kappa_k$"
and is explicitly given 
in essentially the same way as in (\ref{ode_phik.eq}), (\ref{def.Ak}), and (\ref{def.Bk}). 
For the reader's convenience, 
we will give explicit formal expressions of $\phi^{\kappa_k}$ for $k =0, 1, 2, 3$ when $1/3 < H<1/2$.
\begin{eqnarray}
d \phi^{0}_t &=&  \sigma( \phi^{0}_t) dh_t
\qquad
\mbox{ with }
\quad
 \phi^{0}_0 = a, 
  \label{ode_fr0.eq}
\\
d \phi^1_t - \nabla \sigma (  \phi^0_t)  \la \phi^1_t, dh_t \ra
&=&
\sigma (\phi^0_t) dx_t
\quad
\mbox{ with }
\quad
 \phi^{1}_0 = 0,  
 \label{ode_fr1.eq}
 \\
 d \phi^2_t - \nabla \sigma (  \phi^0_t)  \la \phi^2_t, dh_t \ra
&=&
\nabla\sigma (\phi^0_t)\la \phi^1_t, dx_t \ra 
\nn
\\
&&
\quad
+
\frac12 \nabla^2 \sigma ( \phi^0_t) \la \phi^1_t, \phi^1_t,  dh_t\ra
\quad
\mbox{ with }
\quad
 \phi^{2}_0 = 0,
 \label{ode_fr2.eq} 
 \\
 d \phi^{1/H}_t - \nabla \sigma (  \phi^0_t)  \la \phi^{1/H}_t, dh_t \ra
&=&
b (\phi^0_t) dt
\quad
\mbox{ with }
\quad
 \phi^{1/H}_0 = 0.  
 \label{ode_fr1/H.eq}
 \end{eqnarray}

Proposition \ref{pr.map.rp} holds still true  with a slight modification. 
Namely, if $1/p +1/q >1$, the map
\begin{eqnarray*}
\lefteqn{
G\Omega_p ({\cal V} ) \times C_0^{q-var}([0,1], {\cal V} ) 
\ni
({\bf x},h) 
}
\\
&\mapsto&
 ({\bf x}, {\bf h}, \bm{\lambda} ,\tilde{{\bf y}}^{\ve}, \bm{\phi}^0,
  \bm{\phi}^{\kappa_{1}}, \ldots, \bm{\phi}^{\kappa_{k}}, {\bf r}_{\ve}^{\kappa_{k+1}} )
\in
G\Omega_p ({\cal V}^{ \oplus 2} \oplus {\mathbb R} \oplus {\cal W}^{ \oplus k+3}).
\nn
\end{eqnarray*}
is locally Lipschitz continuous for any $k$.
%


The deterministic estimates for  terms in the expansion
(Proposition \ref{pr.map.rp2})
can easily be modified as follows
(This proposition was already used in \cite{inaaop});

\begin{pr}\label{pr.map.rp2.5}
Assume $1/3<1/p <H<1/2$ and $1/p +1/q >1$.
Consider RDE (\ref{rde_fr.eq}) and keep the same notations as above.
Then, the following {\bf (i)} and {\bf (ii)} hold.
\\
\noindent
{\bf (i)}~
For any $\rho >0$ and $k=1,2,\ldots$, there exists a positive constants $C=C(\rho, k)$ 
which satisfies that
\[
\| (\bm{\phi}^{\kappa_{k}} )^1 \|_{p -var} \le C (1+ \overline{\omega_{ {\bf x}}}^{1/p})^{\kappa_{k}}.
\]
for any ${\bf x} \in G\Omega_p ({\cal V} )$
and $h \in C_0^{q-var}([0,1], {\cal V} ) $ with $\| h\|_{q-var} \le \rho$.
\\
\noindent
{\bf (ii)}~
For any $\rho_1, \rho_2 >0$ and $k=1,2,\ldots$,
there exists a positive constants $\tilde{C}= \tilde{C}(\rho_1, \rho_2, k)$, which is independent of $\ve$
and satisfies that
\[
\| ({\bf r}^{\kappa_{k+1}}_{\ve})^1 \|_{p -var} 
\le \tilde{C} (\ve + \ve \overline{\omega_{  {\bf x}}}^{1/p})^{\kappa_{k+1}}
\]
for any ${\bf x} \in G\Omega_p ({\cal V} ) $ with 
$\overline{\omega_{ \ve {\bf x}}}^{1/p}   =  \ve \overline{\omega_{ {\bf x}}}^{1/p}   \le \rho_1$
and any $h \in C_0^{q-var}([0,1], {\cal V} )$ with $\|h\|_{q -var} \le \rho_2$.
\end{pr}

The moment estimates for  terms in the expansion
(Theorem \ref{thm.mom_r})
can be modified in the following way.
This can be shown in essentially the same way as in Theorem \ref{thm.mom_r}.

\begin{tm}\label{thm.mom_r_fr}
We consider RDE (\ref{rde_fr.eq}).
Assume $1/3<1/p <H<1/2$ and $1/p + 1/q >1$
and
let $h \in C_0^{q -var} ([0,1], {\cal V} )$.
Assume that 
${\bf x}$ be a $G\Omega_p ({\cal V})$-valued random variable 
such that 
{\rm (i)}~ $\overline{\omega_{{\bf x}}}= \omega_{{\bf x}} (0,1) \in \cap_{1 \le r <\infty} L^r$ 
and 
{\rm (ii)}~ 
$\exp ( N_{\alpha} ({\bf x}) )
 \in  \cap_{1 \le r <\infty} L^r$ for any $\alpha >0$ .

Then, for any ${\bf x}$, $h$ and $k \in {\mathbb N}$, 
there exist control functions $\eta_k = \eta_{k,{\bf x}, h }$
such that the following {\rm (i)}-- {\rm (iii)} hold:
\\
\noindent
{\rm (i)}~$\eta_{k}$ are non-decreasing in $k$, i.e., 
$\eta_{k,{\bf x}, h}(s,t) \le \eta_{k+1,{\bf x}, h} (s,t)$ 
for all $k, {\bf x}, h, (s,t)$.
\\
\noindent
{\rm (ii)}~$\overline{\eta_{k,{\bf x}, h} } \in  \cap_{1 \le r <\infty} L^r$ 
for all $k, h$.
\\
\noindent
{\rm (iii)}~
For all $\ve \in (0,1]$, $k\in {\mathbb N}$, $h$, ${\bf x}$, 
and $0 \le s \le t \le 1$, $j =1,2$,
we have 
\[
\bigl|  
\bigl(
{\bf x}, {\bf h}, {\bf y}^{\ve}, \bm{\phi}^0, \bm{\phi}^{\kappa_1},\ldots, 
\bm{\phi}^{\kappa_k}, \ve^{- \kappa_{k+1}} {\bf r}_{\ve}^{ \kappa_{k+1} }
\bigr)^j_{s,t}
\bigr| \le \eta_{k,{\bf x}, h} (s,t)^{j/p}.
\]

In particular, for all $k\in {\mathbb N}$ and $h$, 
$\| (\bm{\phi}^{\kappa_{k}})^1\|_{p -var} \in \cap_{1 \le r <\infty} L^r$
and 
$ \| ({\bf r}_{\ve}^{\kappa_{k+1}})^1 \|_{p -var} = O (\ve^{\kappa_{k+1}})$ in $L^r$ 
for any $1 \le r <\infty$.
\end{tm}

\begin{re}\label{re.end_sec2}
{\rm (i)}~ This section (Section \ref{sec_two}) may look a little bit lengthy.
But, we will only use Proposition \ref{pr.map.rp2.5} and Theorem \ref{thm.mom_r_fr}  in later sections.
\\
{\rm (ii)}~The author guesses that 
the results in this section naturally extends to the case of $[p] \ge 3$.
But, computation may be hard and it has not been confirmed yet.
\end{re}

%
%
\section{Malliavin Calculus for solution of RDE driven by fBM}\label{sec.malli}

In this section 
we study the solution of a (scaled) RDE driven by 
fractional Brownian motion with $H \in (1/3,1/2]$ via Malliavin calculus.
It was already done by Hairer and Pillai \cite{hp}
(and Cass, Hairer, Litterer, and Tindel \cite{chlt}).
In this section we basically follow their arguments,
but in our case we need to check dependency on the small parameter $\ve \in (0,1]$.


To keep our argument concise, 
we do not explain much about Malliavin calculus here.
The reader should refer to  well-known textbooks such as 
Nualart \cite{nu} and Shigekawa \cite{sh}.
In this paper we use Watanabe distribution theory and 
asymptotic theorems for them,
which can be found in \cite{wa} or Section V-9, \cite{iwbk}.
(The results in \cite{wa, iwbk} are formulated on the classical Wiener space,
but they are still true on an abstract Wiener space.)
One thing different from is \cite{wa, iwbk} that
the index sets  of asymptotic expansions may not be ${\mathbb N}=\{0,1,2,\ldots\}$
in this paper.
So, we need to slightly modify  these asymptotic theorems.
However, we skip details here 
since a summary was already given in the author's previous work \cite{ina1}.

In this paper, we use the following notations.
$D$ stands for the ${\cal H}$-derivative.
Sobolev space of the integral index $r \in (1,\infty)$
and 
the differential index $s \in {\mathbb R}$ 
is denoted by 
${\bf D}_{r,s}$.
As in \cite{wa, iwbk}, we set 
${\bf D}_{\infty}= \cap_{k=1}^{\infty} \cap_{1 <r<\infty} {\bf D}_{r,k}$,
${\bf D}_{-\infty}= \cup_{k=1}^{\infty} \cup_{1 <r<\infty} {\bf D}_{r, -k}$.
Moreover, we also use
$\tilde{\bf D}_{\infty}= \cap_{k=1}^{\infty} \cup_{1 <r<\infty} {\bf D}_{r,k}$
and 
$\tilde{\bf D}_{-\infty}= \cup_{k=1}^{\infty} \cap_{1 <r<\infty} {\bf D}_{r, -k}$
in Watanabe distribution theory.
The Sobolev space of vector-valued Wiener functionals is denoted by
${\bf D}_{r,s}({\cal K})$, etc., 
where ${\cal K}$ is a real separable Hilbert space.


Let $1/3<H \le 1/2$ and choose $p$ so that $1/3< 1/p <H$.
The $d$-dimensional fBm  $(w_t)_{0 \le t \le 1}$ with Hurst parameter $H$
admits a natural rough path lift ${\bf w}$ as a random rough path  
that takes values in $G\Omega_p  ({\mathbb R}^d)$. 
We denote by ${\cal H}= {\cal H}^H$ the Cameron-Martin space
associated with $d$-dimensional fBm with $H \in (1/3, 1/2]$.
Throughout this section $\gamma \in {\cal H}$ is arbitrary, but fixed.
By Friz-Victoir \cite{fv06}, 
there is a continuous embedding
\begin{equation}\label{emb_fv.eq}
{\cal H}^H \hookrightarrow W^{1/q, 2}( [0,1], {\mathbb R}^d) \hookrightarrow  C_0^{q- var}([0,1],{\mathbb R}^d)
\end{equation}
for any $q \in ( (H+1/2)^{-1}, 2)$. 
(In a very recent preprint \cite{fggr}, 
the above embedding is shown to still hold for $q =(H+1/2)^{-1}$.)
The Banach space in the middle is the fractional Sobolev (i.e., Besov) space with
the differential index  $1/q$ and the integral index $2$.
Note that 
if $p$ and $q$ are sufficiently close to $1/H$ and  $(H+1/2)^{-1}$, respectively,
then $1/p +1/q >1$,
which makes Young integration/translation/paring possible.

Let us make a remark on H\"older regularity of the above RDE.
It is well-known that ${\bf w}$ is actually an $\alpha$-H\"older geometric rough path a.s.,
where we set $\alpha :=1/p$.
At first, it is not so obvious whether $\tau_{\gamma} (\ve {\bf w})$ is an
$\alpha$-H\"older geometric rough path, 
even though ${\cal H}^H \hookrightarrow C_0^{\alpha- hld}({\mathbb R}^d)$.
It was shown to be true
in Friz and Victoir \cite{fv06} and Excercise 9.37, p. 211, \cite{fvbk}.
%
%
%
Using (\ref{emb_fv.eq}) with $1/q =\alpha +1/2$, they showed that
\begin{eqnarray*}
\|\tau_{\gamma} ( {\bf x})^1\|_{\alpha -hld} &\le& \mbox{const.} \times  
( \|{\bf x}^1\|_{\alpha -hld} + \|\gamma\|_{{\cal H}} ),
\\
\| \tau_{\gamma} ({\bf x})^2 \|_{2\alpha -hld} &\le& \mbox{const.} \times
(  \|{\bf x}^2\|_{2\alpha -hld} 
+  \|{\bf x}^1\|_{\alpha -hld}\|\gamma\|_{{\cal H}} + \|\gamma\|_{{\cal H}}^2 )
\end{eqnarray*}
for any $\gamma \in {\cal H}$ and ${\bf x} \in G\Omega_{\alpha}^H  ({\mathbb R}^d)$.
These imply that the driving signal $(\tau_{\gamma} (\ve {\bf w})  , \ve^{1/H} \lambda  )$
 of RDE (\ref{rde_mal.eq}) 
is actually a $\alpha$-H\"older geometric rough path a.s.
Consequently,  so is $\tilde{{\bf y}}^{\ve}$.


As before 
$\sigma : {\mathbb R}^n \to {\rm Mat}(n,d)$ 
and $b : {\mathbb R}^n \to {\mathbb R}^n$ be $C_b^{\infty}$.
For notational convenience, 
we will sometimes denote by $V_i :{\mathbb R}^n \to {\mathbb R}^n$ the $i$th column 
vector field of $\sigma$ ($1 \le i \le d$), i.e., $\sigma =[V_1; \cdots; V_d]$. 
In a similar way we will write $V_0 =b$.

We consider the following RDE for $\ve \in (0,1]$ and $a \in {\mathbb R}^n$;
\begin{eqnarray}\label{rde_mal.eq}
d \tilde{ y}^{\ve}_t 
&=& 
\sigma (\tilde{y}^{\ve}_t) ( \ve d w_t + d \gamma_t)
+ \ve^{1/H} b (\tilde{y}^{\ve}_t) dt
\qquad
\mbox{ with }
\quad
\tilde{ y}^{\ve}_0 = a \in {\mathbb R}^n.
\end{eqnarray}
We write $\tilde{y}^{\ve}_t = a + (\tilde{{\bf y}}^{\ve})^1_{0,t}$
and study this process.
When $\gamma =0$, we write $\tilde{{\bf y}}^{\ve} = {\bf y}^{\ve}$.
When $\gamma =0$ and $\ve =1$, 
we write $\tilde{{\bf y}}^{\ve} = {\bf y}$.
If $\Phi$ denotes the Lyons-It\^o map that corresponds to $[\sigma, b]$
and $a$, then $\tilde{y}^{\ve} = \Phi ( ( \tau_{\gamma} (\ve {\bf w})  , \ve^{1/H} \lambda ))$.
Here, 
{\rm (i)}~$\tau_{\gamma} (\ve {\bf w}) $ 
denotes the Young translation of $\ve {\bf w}$ by $\gamma$
and {\rm (ii)}~
$( \tau_{\gamma} (\ve {\bf w})  , \ve^{1/H} \lambda )$ denotes the Young pairing
of $\tau_{\gamma} (\ve {\bf w})$
and the one-dimensional path $\ve^{1/H} \lambda_t =\ve^{1/H}t$.
Using $V_i$'s we can rewrite RDE (\ref{rde_mal.eq}) as follows:
\begin{eqnarray}\label{rde_mal2.eq}
d \tilde{ y}^{\ve}_t &=& \sum_{i=1}^d V_i (\tilde{y}^{\ve}_t) 
( \ve d w_t^{i} + d \gamma_t^{i})
+ \ve^{1/H} V_0 (\tilde{y}^{\ve}_t) dt
\qquad
\mbox{ with }
\quad
\tilde{y}^{\ve}_0 = a \in {\mathbb R}^n.
\end{eqnarray}
%
%
%
Note that $(y^{\ve}_t)_{0 \le t \le 1}$ and $(y_{\ve^{1/H} t})_{0 \le t \le 1}$
have the same law. 
(See Inahama \cite{inaaop} for a proof).




In Hairer and Pillai \cite{hp},  they proved the following:  
{\rm (i)}~$y_t \in {\bf D}_{\infty} ({\mathbb R}^n)$ for any $t>0$, i.e., 
$D^m y_t$ exists and in $\cap_{1<r<\infty} L^r$  for any $m=0,1,2,\ldots$.
{\rm (ii)}~Under H\"ormander's hypoellipticity condition on vector fields
$\{V_1, \ldots, V_d, V_0\}$ at the starting point $a$, 
Malliavin covariance matrix of $y_t$ 
is non-degenerate in the sense of Malliavin for any $t>0$, i.e., 
\[
  \det \Bigl[\{  \la Dy^{(i)}_t ,Dy^{(j)}_t\ra_{{\cal H}}  \}_{i,j =1}^n    \Bigr]^{-1} 
  \in   \cap_{1<r<\infty} L^r,
\]
where $y^{(i)}_t$ denoted the $i$th component of $y_t$.

It is almost obvious that $\tilde{y}^{\ve}_1$ also satisfies 
{\rm (i)} and {\rm (ii)} above for each fixed $\ve$.
In this paper, however, 
we need to check dependency on $\ve \in (0,1]$ as it varies.
The precise statements are given in the following two propositions.
We will prove them later by slightly modifying the proofs in \cite{hp, chlt, ina3}.

\begin{pr}\label{pr.d_infty}
Assume $\sigma$ and $b$ are $C_b^{\infty}$ and let $\gamma \in {\cal H}$ be arbitrary but fixed. 
Then, for any $m=0,1,2,\ldots$ and  $r \in (1, \infty)$, 
there exists a positive constant $c=c_{m,r}$ such that
\[
{\mathbb E} \bigl[ \| D^m \tilde{y}^{\ve}_1 \|_{ {\cal H}^{\otimes m} }^r  \bigr]^{1/r}  \le c \ve^m.
\]
\end{pr}

\Proof
In \cite{ina3} 
the author proved ${\bf D}_{\infty}$-property of solutions of RDEs 
driven by Gaussian rough path ${\bf w}$ including fBm with $H >1/4$.
The proof is so flexible that we can replace 
${\bf w}$ by $\tau_{\gamma} (\ve {\bf w}) = \ve {\bf w} +\gamma$.
If we keep track of $\ve$-dependency in that argument, 
then we can easily see that
$D^m \tilde{y}^{\ve}$ is $O(\ve^m)$ as $\ve \searrow 0$ for any $m \in {\mathbb N}$.
In that proof,
the uniform estimate of Jacobian process and its inverse 
plays a crucial role.
\QED

\begin{pr}\label{pr.nondeg}
In addition to the assumption of Proposition \ref{pr.d_infty}, 
we assume the ellipticity assumption {\bf (A1)}.
Then, $(\tilde{y}^{\ve}_1 -a)/\ve$ is uniformly non-degenerate in the sense of Malliavin,
that is,
\[
\sup_{0<\ve \le 1}
 {\mathbb E} \Bigl[
  \det \bigl[\{  \la D \Bigl( \frac{\tilde{y}^{\ve, (i)}_1 -a}{\ve} \Bigr) , 
   D \Bigl( \frac{\tilde{y}^{\ve, (j)}_1 -a}{\ve} \Bigr) 
      \ra_{{\cal H}}  \}_{i,j =1}^n    \bigr]^{-r} 
  \Bigr]
  <\infty
\]
for any $r \in (1, \infty)$.
\end{pr}

\Proof
Note that the special case 
"$\gamma =0$ and $b \equiv 0$ and uniformly elliptic coefficients"
was already shown in \cite{bnot, boz1}, etc.
Since this proposition can be shown in a similar way,  
we omit the proof.
(However, we note that uniform non-degeneracy of $(\tilde{y}^{\ve}_1 -a)/\ve$ 
for the shifted RDE becomes quite complicated 
under a H\"ormander-type condition 
instead of {\bf (A1)}.)

\QED


Consider the asymptotic expansion of $\tilde{y}^{\ve}$ as in (\ref{def.r.2}).
We have already seen that this expansion holds true 
both in the deterministic sense and the $L^r$-sense.
Moreover, evaluated at time $t=1$, 
it also holds true in ${\bf D}_{\infty}$-sense.

\begin{pr}\label{pr.asyDinf}
We keep the same assumptions as in Proposition \ref{pr.d_infty}.
Then, we have the following asymptotic expansion as $\ve \searrow 0$:
\[
\tilde{y}^{\ve}_1 \sim
 \phi^0_1  + \ve^{\kappa_1} \phi^{\kappa_1}_1 +\cdots 
+ \ve^{\kappa_k} \phi^{\kappa_k}_1 + \cdots
\qquad
\qquad
\mbox{ in ${\bf D}_{\infty} ({\mathbb R}^n)$.}
\]
This means that for each $k$, 
{\rm (i)} $\phi^{\kappa_k}_1 \in {\bf D}_{\infty} ({\mathbb R}^n)$ 
and  {\rm (ii)}
${\bf D}_{r,s}$-norm of  
$r_{\ve, 1}^{ \kappa_{k+1}}$ is $O (\ve^{\kappa_{k+1}} )$ for any $r \in (1, \infty)$
and $s \ge 0$.
\end{pr}

\Proof 
By the way it is constructed, $\phi^{\kappa_k}_1$ is an element of 
inhomogeneous Wiener chaos of order at most $[\kappa_k]$.
Hence, $\phi^{\kappa_k}_1 \in {\bf D}_{\infty} ({\mathbb R}^n)$
and  $D^{[\kappa_k] +1} \phi^{\kappa_k}_1 =0$.
Next we estimate Sobolev norms of the remainder terms. 
We see from the stronger form of Meyer's equivalence that,
for any integer $s \ge [\kappa_k] +1$ and any $r \in (1, \infty)$,
there exists $C=C_{r,s}$ such that
\[
\| r_{\ve, 1}^{  \kappa_{k+1} }\|_{{\bf D}_{r,s}} 
\le C ( \| r_{\ve, 1}^{\kappa_{k+1} }\|_{L^r} +  \|D^s r_{\ve, 1}^{ \kappa_{k+1}}\|_{L^r})
=
C ( \| r_{\ve, 1}^{\kappa_{k+1} }\|_{L^r} +  \|D^s   \tilde{y}^{\ve}_1\|_{L^r})
\]
holds.
By Theorem \ref{thm.mom_r_fr} and Proposition \ref{pr.d_infty},
the right hand side is $O(\ve^{\kappa_{k+1}}) + O(\ve^s)= O(\ve^{\kappa_{k+1}})$.
Thus, we have the desired estimate for such $(r,s)$.
Since ${\bf D}_{r,s}$-norm  is increasing in $s$, the proof is done.
\QED


Now we state and prove  on-diagonal short time asymptotics of 
$p_t (a,a) ={\mathbb E} [ \delta_a ( y_t)]$.
Compared to the off-diagonal case, this is not so difficult.
From Propositions \ref{pr.nondeg}, and \ref{pr.asyDinf},
and Watanabe's asymptotic theory for generalized Wiener functionals (i.e., Watanabe distributions), 
we can obtain the following theorem.
\begin{tm}\label{tm.ondiag}
Assume the ellipticity assumption {\bf (A1)}.
Then, the diagonal of the kernel $p(t, a, a)$ admits the following 
asymptotics as $t \searrow 0$;
\[
p(t, a, a) \sim
 \frac{1}{t^{n H}} \bigl( c_0 + c_{\nu_1} t^{\nu_1 H} +  c_{\nu_2} t^{\nu_2 H} +  \cdots \bigr)
\]
for certain real constants $c_0, c_{\nu_1}, c_{\nu_2}, \ldots$.
Here, $\{ 0=\nu_0 <\nu_1<\nu_2 <\cdots \}$ are all the elements of 
$\Lambda_3$ in increasing order.
Moreover, $c_0$ is positive.
\end{tm}

\Proof 
In this proof, $\gamma =0$.
From the scaling property, we see that 
$$ 
p( \ve^{1/H}, a, a) =  {\mathbb E}[ \delta_{a}  ( y^{\ve}_1 (a)) ]
=
 {\mathbb E}[ \delta_{0}    \bigl( \ve \frac{ y^{\ve}_1 (a) -a}{\ve}  \bigr) ]
=
\ve^{-n} {\mathbb E}[ \delta_{0}    \bigl(  \frac{ y^{\ve}_1 (a) -a}{\ve}  \bigr) ].
$$

By Proposition \ref{pr.nondeg}, $( y^{\ve}_1 (a) -a)/\ve$
is uniformly non-degenerate.
It admits asymptotic expansion in ${\bf D}_{\infty} ({\bf R}^n)$ as in Proposition \ref{pr.asyDinf}
with the index set for the exponents being $\Lambda_2$.
Then, by (a slight generalization of)
Theorem 9.4, p. 387, Ikeda and Watanabe \cite{iwbk}, the following asymptotic expansion holds in 
$\tilde{{\bf D}}_{- \infty}$ as $\ve \searrow 0$;
\[
\delta_{0}    \bigl(  \frac{ y^{\ve}_1 (a) -a}{\ve}  \bigr)
\sim
\phi_0 + \ve^{\nu_1} \phi_{\nu_1} + \ve^{\nu_2} \phi_{\nu_2} +\cdots
\qquad
\mbox{as $\ve \searrow 0$.}
\]
Formally, this is a composition of Taylor expansion of $\delta_0 (\,\cdot\,)$ 
and the asymptotic expansion of $( y^{\ve}_1 (a) -a )/\ve$.
Hence, the new index set is ${\mathbb N} \la \Lambda_2\ra = \Lambda_3$.
By taking the generalized expectation and setting $c_{\nu_k} = {\mathbb E}[\phi_{\nu_k} ]$,  we have
\[
p( \ve^{1/H}, a, a) 
\sim
\ve^{-n}
\bigl(
c_0 + 
c_{\nu_1}  \ve^{\nu_1} 
+ c_{\nu_2} \ve^{\nu_2}  +\cdots \bigr)
\qquad
\mbox{as $\ve \searrow 0$.}
\]
Putting $\ve = t^H$, we  prove the asymptotic expansion.
It is straight forward to see that 
$$
c_{0} = {\mathbb E}[\delta_{0} ( \sum_{j=1}^d V_j (a) w_1^j) ]  
= (2\pi)^{-n/2} \{ \det (\sigma (a) \sigma (a)^*) \}^{-1/2} >0,
$$
which completes the proof.
\QED

%
%
%
%
%
%
\section{Off-diagonal short time asymptotics}
\label{sec.pf}

In this section, following Watanabe \cite{wa}, 
 we prove the short time asymptotics of kernel function
$p_t(a,a')$ when $a \neq a'$ and $1/3 <H \le 1/2$.
%
%
Unlike in \cite{wa},  
we can localize around the energy minimizing path in the geometric rough path space
in this paper,
since Lyons-It\^o map is continuous in this setting.
(The case $H >1/2$ was done in \cite{ina1}.
The result in this section can be regarded as 
a rough path version of that in \cite{ina1}.)

\subsection{Localization around energy minimizing path}
\label{subsec.ldp}

Let $G\Omega^B_{\alpha, m} ({\mathbb R}^d)$ be the geometric rough path space with 
$(\alpha, m)$-Besov norm for $\alpha \in (1/3, 1/2]$ and $m>1$ with $\alpha -1/m >1/3$.
Explicitly, the norms are given by
\[
\| {\bf x}^i \|_{i\al, m/i -B} :=  \Bigl(  \iint_{0 \le s < t \le 1}  
\frac{ | {\bf x}^i_{s,t} |^{m/i} }{ |t-s|^{1 +m \al} }dsdt \Bigr)^{i/m}
\qquad
\qquad
(i=1,2).
\]
We have the following continuous embeddings
\begin{equation}\label{go_emb.eq}
 G\Omega^H_{\beta} ({\mathbb R}^d)
\hookrightarrow  G\Omega^B_{\alpha', m} ({\mathbb R}^d) 
\hookrightarrow  G\Omega^H_{\al} ({\mathbb R}^d)
\hookrightarrow  G\Omega_{p} ({\mathbb R}^d)
\end{equation}
if $1/3 <1/p =\al < \al' -1/m <\al' < \beta \le 1/2$
(see Appendix A2, Friz and Victoir \cite{fvbk}).

Next, we introduce a measure.
Let $\mu =\mu^H$ be the law of the fractional Brownian motion 
with Hurst parameter $H \in (1/3, 1/2]$.
This is a probability measure on ${\cal W}= \overline{{\cal H}}$,
which is the closure of ${\cal H}= {\cal H}^H$ in 
$C_0^{p -var} ([0,1], {\mathbb R}^d)$.
Then, the triple 
$({\cal W}, {\cal H}, \mu)$ is an abstract Wiener space.

For any $\beta \in (1/3, H)$,
fBm $(w_t)$ admits a natural lift a.s. via dyadic piecewise linear approximation  
and the lift ${\bf w}$ is a random variable taking values in $G\Omega^H_{\beta} ({\mathbb R}^d)$.
Note that the lift of Cameron-Martin space ${\cal H}$
is contained in $G\Omega^H_{\beta} ({\mathbb R}^d)$.
Moreover, as $\ve \searrow 0$,
Schilder-type large deviation holds for the laws of $\ve {\bf w}$,
which will be denoted by $\nu_{\ve} = \nu^H_{\ve}$.
(See Friz and Victoir \cite{fv07}).
Because of Besov-H\"older embedding mentioned above, 
these properties also hold with respect to $(\al', m)$-Besov topology if $\al' <H$.
As usual,
the good rate function ${\cal I}$ is given as follows:
${\cal I} ({\bf x})= \|h\|_{\cal H}^2 /2$ if ${\bf x}$ is the lift 
of some $h \in {\cal H}$
and 
${\cal I} ({\bf x})= \infty$ if otherwise.

Let us clarify the conditions on various indices here.
From now on, these will be assumed unless otherwise stated.
First, for given $H \in (1/3, 1/2]$,
 we choose $p :=1/\al \in (1/H,3)$ and $q\in ( (H +1/2)^{-1},2)$ so that  
$1/p +1/q >1$ holds.
Then, 
we choose $\al' \in (\al, H)$ and $m \in {\mathbb N}$
such  that
$(\al' -\al) \vee (H -\al') > 1/(4m)$
and consider $G\Omega^B_{\alpha', 4m} ({\mathbb R}^d)$.
(Heuristically, $m$ is a very large integer.)

Since $m$ is an integer, 
$w \mapsto \| {\bf w}^i - {\bf h}^i \|_{i\al, 4m/i -B}^{4m/i}$ 
is 
 ${\bf D}_{\infty}$ in the sense of Malliavin calculus for $i=1,2$,
 where ${\bf h}$ is a fixed element.
 Actually, it is an element of an inhomogeneous Wiener chaos.
Due to this fact, the localization is allowed even in the framework of 
 Watanabe distribution theory.
This is the reason why we use this Besov-type norm on the geometric rough path space.

The Young translation $\tau_{\gamma}$ works on $G\Omega^B_{\alpha', 4m} ({\mathbb R}^d) $
for any $\gamma \in {\cal H}$.
The proof is just a slight modification of the H\"older case.
\begin{lm}\label{lm.YtrBes}
Let $H, \al', m$ be as above.
Then, for any $\gamma \in {\cal H}$, 
the Young translation $\tau_{\gamma}$ is a continuos map from 
 $G\Omega^B_{\alpha', 4m} ({\mathbb R}^d) $ to itself.
\end{lm}

\Proof 
Generally, we have the following basic result for Young integrals.
Let  $p' , q' >0$  with $1/p' +1/q' >1$.
Then, there is a constant $C>0$ which depends only on $p' , q'$
such that
\[
\Bigl|
\int_s^t
(x_u - x_s) \otimes dy_u
\Bigr| 
\le 
C \|x\|_{p' -var; [s,t]} \cdot \|y\|_{q' -var; [s,t]}.
\]
for any $[s,t] \subset [0.1]$.

Now we prove the lemma.
We have
\begin{align}
\tau_{\gamma} ({\bf x})^1_{s,t} &=  {\bf x}^1_{s,t} + \bm{\gamma}^1_{s,t}
\nn
\\
\tau_{\gamma} ({\bf x})^2_{s,t} &=  {\bf x}^2_{s,t} + \bm{\gamma}^2_{s,t}
+
\int_s^t
{\bf x}^1_{s,u} \otimes d\gamma_u 
+
\int_s^t
\bm{\gamma}^1_{s,u} \otimes d x_u 
\label{tau.hyoki.eq}
\end{align}
Here, the second, the third, and the fourth terms on the right hand side of (\ref{tau.hyoki.eq})
are Young integrals.
As usual we set $x_t = {\bf x}^1_{0,t}$.
By Besov-H\"older embedding theorem, $x$ is $\al' -1/(4m)$ H\"older continuous.
Moreover, 
there is a constant $c$ such that
 $$
 \|\gamma\|_{q -var; [s,t]} \le c \| \gamma \|_{W^{1/q,2}} \cdot  (t-s)^{  \frac{1}{q} - \frac12 }
 \le 
 c \| \gamma \|_{{\cal H}} \cdot  (t-s)^{ \frac{1}{q} - \frac12 }
 \qquad
 (\gamma \in {\cal H}, \quad \frac{1}{q} < H +\frac12).
$$
(See p. 211, \cite{fvbk}. The constant $c>0$ may vary from line to line.)
Therefore, $\bm{\gamma}^2$ is of finite $2( 1/q -1/2)$ H\"older norm.
The third and 
the fourth terms are of finite $( 1/q -1/2) + (\al' -1/(4m))$ H\"older norm.
Since $H -\al'<  1/(4m)$ and we may choose $q$  so that
$ 1/q -1/2$ can be arbitrarily close to $H$,
these three terms are actually of finite $(2 \al' +\delta)$  H\"older norm
for some $\delta >0$
and hence are of finite $(2\al', 2m)$ Besov norm.
Thus, we have shown that 
$\tau_{\gamma}$ maps  $G\Omega^B_{\alpha', 4m} ({\mathbb R}^d) $ to itself.

We can show continuity of $\tau_{\gamma}$ 
by estimating the difference $|\tau_{\gamma} ({\bf x})^i_{s,t} - \tau_{\gamma} (\tilde{\bf x})^i_{s,t}|$
for $i=1,2$ in essentially the same way. 
So, we omit details.
 \QED


For $\gamma \in {\cal H} \subset C_0^{q -var} ([0,1], {\mathbb R}^d)$,
let $\phi^0= \phi^0(\gamma)$ be a unique solution of  (\ref{ode_fr0.eq})
in the $q$-variational Young sense, 
which starts at $a \in {\mathbb R}^n$.
Set,  for $a \neq a'$, 
\[
K_a^{a'} = \{ \gamma \in {\cal H} ~|~  \phi^0_1(\gamma) =a'\}.
\]
This is a closed set in ${\cal H}$.
We only consider the case that $K_a^{a'} $ is not empty.
For example, if {\bf (A1)} is satisfied for any $a$, then $K_a^{a'} $ is not empty for any $a'$.
From the Schilder-type large deviation theory, we see that
$\inf\{ \|\gamma\|_{\cal H} ~|~ \gamma \in  K_a^{a'}\} 
= \min\{ \|\gamma\|_{\cal H} ~|~ \gamma \in  K_a^{a'}\}$.

We continue to assume {\bf (A1)}.
Now we introduce another assumption;
\vspace{3mm}
\\
{\bf (A2):} $\bar{\gamma} \in K_a^{a'}$ which minimizes ${\cal H}$-norm exists uniquely.
\vspace{3mm}
\\
In the sequel, $\bar{\gamma}$ denotes the minimizer in Assumption {\bf (A2)}
and we use the results of the previous section for this $\bar{\gamma}$.

Note that {\rm (i)}~ the mapping $\gamma \in {\cal H} \hookrightarrow
 C_0^{q -var} ([0,1], {\mathbb R}^d)\mapsto  \phi_1^0 (\gamma) \in {\mathbb R}^n$
is Fr\'echet differentiable
and 
{\rm (ii)}~ its Jacobian is a surjective linear mapping from ${\cal H}$ to ${\mathbb R}^n$ at any $\gamma$,
because there exists a positive constant $c=c(\gamma)$ such that 
\begin{equation}\label{detuki.eq}
\Bigl(  \la D  \phi_1^{0,i } (\gamma),   D  \phi_1^{0, j} (\gamma) \ra_{{\cal H}^*}  \Bigr)_{ 1 \le i,j \le n}  \ge  c \cdot {\rm Id}_n.
\end{equation}
This can be shown in the same way as in the proof of 
non-degeneracy of $y_t$ under ellipticity assumption.
(Actually, it is easier since $\gamma$ is non-random and fixed here.)

Therefore, by the Lagrange multiplier method,  there exists 
$\bar{\nu} =(\bar{\nu}_1, \ldots, \bar{\nu}_n) \in {\mathbb R}^n$ uniquely such that
the map
\begin{equation}\label{Lmul1.eq}
{\cal H} \times {\mathbb R}^n \ni (\gamma, \nu) \mapsto   
\frac12   \|\gamma\|_{\cal H}^2 - \la  \nu, \phi^0_1(\gamma) -a' \ra_{{\mathbb R}^n} \in {\mathbb R}
\end{equation}
attains an extremum at $(\bar{\gamma}, \bar{\nu})$.
By differentiating in the direction of $k \in {\cal H}$, 
we have
\begin{equation}\label{Lmul2.eq}
\la \bar{\gamma} , k \ra_{\cal H} = \la  \bar{\nu},  D_k \phi^0_1( \bar{\gamma} ) \ra_{{\mathbb R}^n}
=
\bigl\la  \bar{\nu},  
\hat{J}(\bar{\gamma})_1 \int_0^1   \hat{J}(\bar{\gamma})_t^{-1}     \sigma(   \phi^0_t(\bar{\gamma} ) )  dk_t
 \bigr\ra_{{\mathbb R}^n}.
\end{equation}
Here, 
$\hat{J}(\bar{\gamma})^{\pm 1}$ are of finite  $q$-variation
and $\hat{J}(\bar{\gamma})$ satisfies the following ODE in Young sense;
\[
dJ_t = \nabla \sigma (\phi^0_t ( \bar{\gamma})) ) \la J_t  , d \bar{\gamma}_t \ra
\qquad
\mbox{ with $J_0 ={\rm Id}_n$.}
\]
%
Since the integral on the right hand side is of (\ref{Lmul2.eq}) Young integral,
 $\la \bar{\gamma} , \,\cdot\, \ra_{\cal H}$ naturally extends to a continuous linear 
functional on $C_0^{p -var} ([0,1], {\mathbb R}^d)$.

Next, set $\hat\nu_{\ve} = \nu_{\ve} \otimes \delta_{ \ve^{1/H} \lambda }$,
where $\lambda$ is a one-dimensional path defined by $\lambda_t =t$
and $\otimes$ stands for the product of probability measures. 
%
%
%
This measure is supported on 
$G\Omega^B_{\alpha', 4m} ({\mathbb R}^d)  \times {\mathbb R} \la \lambda \ra$.
The Young pairing map
$G\Omega^B_{\alpha', 4m} ({\mathbb R}^d)  \times {\mathbb R} \la \lambda \ra
\to  G\Omega^B_{\alpha', 4m} ({\mathbb R}^{d+1}) $ is continuous.
The law of $\hat\nu_{\ve}$ induced by this map is 
the law of $(\ve {\bf w}, \ve^{1/H} \lambda)$, the Young pairing of $\ve {\bf w}$ and $\ve^{1/H} \lambda$.

Define $\hat{ {\cal I}} ({\bf x} ; l) 
= \|h\|_{{\cal H}}^2/2$ if ${\bf x}$ is the lift of some $h\in {\cal H}$ and $l_t \equiv 0$
and define $\hat{ {\cal I}} (w, l) =\infty$ if otherwise.
Here, $l$ is a one-dimensional path.
We can easily show that
 $\{\hat\nu_{\ve} \}_{\ve >0}$ also satisfies a large deviation principle 
as $\ve \searrow 0$ with a good rate function $\hat{ {\cal I}}$.
We will use  this in Lemma \ref{lm.ldpcut} below to show that
we may localize on a neighborhood of 
the minimizer $\bar{\gamma}$ in order to obtain the asymptotic expansion.

Now we introduce a cut-off function. 
Let $\psi : {\mathbb R} \to [0,1]$ be a smooth function such that 
$\psi (u) =1$ if $|u| \le 1/2$ and $\psi (u) =0$ if $|u| \ge 1$.
For each  $\eta >0$ and $\ve >0$, we set 
\[
\chi_{\eta} (\ve, w) 
= 
\prod_{i=1}^2
\psi \Bigl(  
\frac{1}{\eta^{4m}}     
\| \tau_{- \bar{\gamma}} (\ve {\bf w} )^i \|_{ i\alpha', 4m/i -B}^{4m/i}
\Bigr).
\]
Here, $\tau_{ - \bar{\gamma}}$ is the Young translation by $- \bar{\gamma}$.
It is a continuous map from $G\Omega^H_{\beta} ({\mathbb R}^d)$ to itself.
So, the right hand side is defined for almost all $w \in {\cal W}$.
Shifting by $\bar{\gamma}/\ve$, we have
$$\chi_{\eta} (\ve, w +\frac{\bar{\gamma}}{\ve}) 
=
\prod_{i=1}^2
\psi \Bigl(  
\frac{\ve^{4m}}{\eta^{4m}}     \|  {\bf w}^i \|_{i\alpha', 4m/i -B}^{4m/i}
\Bigr).
$$
This is a ${\bf D}_{\infty}$-functional.
Moreover, from Taylor expansion for $\psi$,
the following asymptotics holds; for any $\eta >0$ and any $M \in {\mathbb N}$, 
\begin{equation}\label{chi_asy.eq}
\chi_{\eta} (\ve, w +\frac{\bar{\gamma}}{\ve}) = 1+O (\ve^M)
\qquad
\mbox{in ${\bf D}_{\infty}$ as $\ve \searrow 0$.}
\end{equation}

Since  $\|  {\bf w}^i \|_{i\alpha', 4m/i -B}^{4m/i}$ is an element of 
an inhomogeneous Wiener chaos of order $4m$, 
so is its Cameron-Martin shift
$\| \tau_{- \bar{\gamma}} (\ve {\bf w} )^i \|_{i\alpha', 4m/i -B}^{4m/i}$.
For any $r \in (0, \infty)$, $L^r$-norm of
this Wiener functional is bounded in $\ve$.
Hence, so is its  ${\bf D}_{r,k}$-norm for any $r, k$.


%
%
The following lemma states that only rough paths sufficiently close to
the lift of the minimizer $\bar{\gamma}$ contribute to the asymptotics.
\begin{lm}\label{lm.ldpcut}
Assume {\bf (A1)} and {\bf (A2)}.
Then, for any $\eta >0$, there exists $c=c_{\eta}>0$ such that 
\[
0 \le {\mathbb E} \bigl[ (1 - \chi_{\eta} (\ve, w) ) \cdot  \delta_{a'} (y^{\ve}_1 )
\bigr]
=
O \Bigl(
\exp \bigl\{  - \frac{   \|  \bar{\gamma}\|_{\cal H}^2 +c }{2 \ve^2}  \bigr\}
\Bigr)
\qquad
\mbox{as $\ve \searrow 0$.}
\]
\end{lm}

\Proof
The proof of this lemma is a bit lengthy and 
quite similar to the proof for the corresponding lemma in \cite{ina1} or \cite{wa},
except that we work on the geometric rough path space.
So, we only give a sketch of proof here.

Set $g(u) = u \vee 0$ for $u \in {\mathbb R}$.
Then, in the sense of distributional derivative, $g^{\prime \prime}  =\delta_0$.
Take a bounded continuous function $C: {\mathbb R}^n \to {\mathbb R}$ such that  
$C(u_1, \ldots, u_n) = g(u_1- a'_1) g(u_2- a'_2)\cdots g(u_n- a'_n)$ 
if $|u -a'| \le 2\eta'$.
Take $\eta' >0$ arbitrarily small.

Then, we have
\begin{align}
0 \le
 {\mathbb E} \bigl[ (1 - \chi_{\eta} (\ve, w) ) \cdot  \delta_{a'} (y^{\ve}_1 )
\bigr]
&=
 {\mathbb E} 
\Bigl[ \{ 1 - \chi_{\eta} (\ve, w) \} 
\psi \Bigl(   
\frac{ | y^{\ve}_1 -  a' |^2  }{\eta'^2 }
\Bigr)
 \cdot  \delta_{a'} (y^{\ve}_1 )
\Bigr]
\nn\\
&=
 {\mathbb E} 
\Bigl[ \{ 1 - \chi_{\eta} (\ve, w) \} 
\psi \Bigl(   
\frac{ | y^{\ve}_1 -  a' |^2  }{\eta'^2 }
\Bigr)
 ( \partial_1^2 \cdots \partial_n^2 C) (y^{\ve}_1 )
\Bigr]
\label{ldp_pf1}
\end{align}
The idea is that by using  
the integration by parts formula for generalized expectations
as in \cite{wa, iwbk},
we reduce the problem to the upper bound of 
the large deviation principle for $\{ \hat{\nu}_{\ve} \}_{\ve >0}$.
(The reason is as follows.
Thanks to the formula, we can remove 
the partial differentiations in 
$( \partial_1^2 \cdots \partial_n^2 C) (y^{\ve}_1 )$
at a certain price.
Then, we have only to treat $C (y^{\ve}_1 )$ 
which is just a bounded function.)
\QED

\subsection{Integrability lemmas}

In this subsection, we prove a few lemmas for integrability of Wiener functionals
of exponential type which will be used in the proof of the short time asymptotic expansion.

Throughout this subsection we assume {\bf (A2)}.
Let $\bar{\gamma}$ be as in {\bf (A2)} and let $\phi^{\kappa_j}$ and 
$r^{\kappa_j +}_{ \ve} =  r^{\kappa_{j +1} }_{ \ve}$ ($j=0,1,2,\ldots$)
be as in (\ref{def.r.2})
 with $\gamma = \bar{\gamma}$. 
%
%
%
First we consider 
\begin{align}
\frac{ r^{\kappa_3}_{ \ve} }{\ve^2}
=
\frac{ r^{2+}_{ \ve} }{\ve^2}
=  
\frac{1}{\ve^2} (\tilde{y}^{\ve} -   \phi^0  -\ve    \phi^1- \ve^2 \phi^2)
%
= \ve^{ \kappa_3 -2 } \phi^{ \kappa_3} + \ve^{ \kappa_4 -2} \phi^{ \kappa_4} +\cdots.
\nn
\end{align}
Recall that $ \kappa_3 = 1/H$ if $H \in (1/3,1/2)$
and $ \kappa_3 = 3$ if $H=1/2$.
When evaluated at time $t=1$, this quantity has a kind of exponential integrability
in the following sense.
(Now that $\bar\gamma$ is fixed, $r^{2+}_{ \ve} ({\bf x})$, $\phi^2 ({\bf x})$, etc. are 
function of ${\bf x}$ alone. 
We will often write $r^{2+}_{ \ve}$, $\phi^2$, etc. for simplicity.)

\begin{lm}\label{lm.integ1}
Assume {\bf (A2)}.
For any $M>0$, there exists $\eta >0$ such that 
\[
\sup_{0 < \ve \le 1} 
{\mathbb E} \bigl[
\exp \bigl(   
M \la \bar\nu,  r^{2+}_{ \ve, 1} \ra /\ve^2
\bigr)
I_{  U_{\eta}   }  (\ve {\bf w}) 
\bigr]  <\infty.
\] 
Here, we set $U_{\eta}= \cap_{i=1,2} 
 \{  \| {\bf x}^i \|_{i\alpha', 4m/i -B}^{1/i}  < \eta \}$ as before.
\end{lm}

\Proof
Let $\omega_{{\bf x}}$ be as in (\ref{nat_control.def}).
Note that $U_1$ is bounded with respect to $p$-variation norm.
%
%
So we may use Proposition \ref{pr.map.rp2} to see that,
for some positive constants $c_1, c_2$,
$$
\| r^{2+}_{ \ve} \|_{p -var} \le c_1 (\ve +  \bar{\omega}_{\ve {\bf x}}^{1/p}  )^{\kappa_3 }
\le 
c_2 (\ve +  \| (\ve  {\bf  x})^1\|_{\alpha' , 4m -B}   
+
   \| (\ve  {\bf  x})^2\|_{2\alpha' , 2m -B}^{1/2}
 )^{\kappa_3 }
\qquad(\ve {\bf x} \in U_{1}).
$$
(In this paragraph we used Besov-H\"older-variation embedding theorem on geometric rough path spaces.
See Proposition A.9,  p. 578, \cite{fvbk} for instance.)
Hence, if  $\ve {\bf x} \in U_{\eta}$ for $0 <\eta \le 1$, then 
\begin{equation} \label{r2+.eq}
\frac{ \| r^{2+}_{ \ve} \|_{p -var} }{\ve^2} 
\le 
c_2 (1 +   \|   {\bf  x}^1\|_{\alpha' , 4m -B}   
+
   \|  {\bf  x}^2\|_{2\alpha' , 2m -B}^{1/2}   )^{2 }
(\ve + 2\eta)^{ \kappa_3 -2 } .
\end{equation}

Recall that  Fernique's theorem holds for fractional Brownian rough path ${\bf w}$
with respect to $\beta$-H\"older topology and hence 
with respect to $(\al' , 4m)$-Besov topology.
It states that for some $\rho >0$ we have 
$$
{\mathbb E} \Bigl[ \exp \Bigl(\rho  (1 +  \|   {\bf  w}^1\|_{\alpha' , 4m -B}   
+
   \|  {\bf  w}^2\|_{2\alpha' , 2m -B}^{1/2} )^2 \Bigr)  \Bigr] <\infty.
$$
(See Friz and Oberhauser \cite{fo} for a proof.)

For given $M$,
take $0< \eta \le 1$ so that $ M |\bar\nu| c_2 (3\eta)^{\kappa_3 -2} \le \rho$.
Then, we have
\[
\sup_{0 < \ve \le \eta} 
{\mathbb E} \bigl[
\exp \bigl(   
M \la \bar\nu,  r^{2+}_{ \ve, 1} \ra /\ve^2
\bigr)
I_{  U_{\eta}   }  (\ve {\bf w}) 
\bigr]  <\infty.
\] 
Note that, 
if $\ve {\bf w} \in U_{\eta}$ and $\eta \le \ve \le 1$, 
then $\| r^{2+}_{ \ve} \|_{p-var} /\ve^2$ is bounded.
(The bound may depend on $\eta$.)
This completes the proof.
\QED


Next we consider 
\begin{align}
\frac{ r^{2}_{ \ve} }{\ve}
=
\frac{ r^{1+}_{ \ve} }{\ve}
=  
\frac{1}{\ve} (\tilde{y}^{\ve} -   \phi^0  -\ve    \phi^1)
= \ve    \phi^{2} +  \ve^{\kappa_3 -1} \phi^{\kappa_3 }
+\cdots.
\nn
\end{align}

\begin{lm}\label{lm.integ2}
Assume {\bf (A2)}.
For any $M>0$, there exists $\eta >0$ such that 
\[
\sup_{0 < \ve \le 1} 
{\mathbb E} \bigl[
\exp \bigl(   
M  \|  r^{2}_{ \ve} \|_{p -var}^2 /\ve^2
\bigr)
I_{  U_{\eta}   }  (\ve {\bf w})
\bigr]  <\infty.
\] 
\end{lm}

\Proof
We can prove the lemma in the same way as in Lemma \ref{lm.integ1} above.
So we only give a sketch of proof.

In this case we have the following inequality instead of (\ref{r2+.eq}):
$$
\frac{ \| r^{2}_{ \ve} \|_{p -var}^2 }{\ve^2} 
\le 
c_2 (1 +   \|   {\bf  x}^1\|_{\alpha' , 4m -B}   
+
   \|  {\bf  x}^2\|_{2\alpha' , 2m -B}^{1/2}   )^{2 }
(\ve + 2\eta)^{ 2 } 
\qquad
(\ve {\bf x} \in U_{\eta}).
$$
The rest is similar. So we omit details.
\QED


From now on we assume {\bf (A1)} and {\bf (A2)}.
In addition, we introduce the following assumption;
\vspace{5mm}
\\
{\bf (A3)':}  \qquad\qquad
${\mathbb E}[  \exp \bigl(  \la\bar\nu ,  \phi^2_1 ({\bf w}) \ra \bigr)  ~|~ \phi^1_1  =0 ]
< \infty.$
\vspace{5mm}

Note that 
$\phi^1_T (w) = \hat{J}_T
 \int_0^T \hat{J}_t^{-1} \sigma ( \phi^0_t) dw_t$.
Here $\phi^0_t = \phi^0_t (\bar\gamma)$, $\hat{J}_t = \hat{J}(\bar\gamma)_t$.
Note that the right hand side is Young integral
and, consequently, is continuous in $w \in {\cal W}$.
%
%
%
%
We regard its $j$th component $\phi^{1, j}_1 \in {\cal W}^* \subset {\cal H}^*$ 
as an element of ${\cal H}$ by Riesz isometry, 
we write ${}^{\sharp} \phi^{1, j}_1 \in {\cal H} \subset {\cal W}$.
We have an orthogonal decomposition
${\cal H} =  
 \ker \phi^1_1  \oplus  ( \ker \phi^1_1 )^{\bot}$.
We denote by $\pi$ the orthogonal projection from ${\cal H}$ onto $\ker \phi^1_1$.
Note that
$ ( \ker \phi^1_1 )^{\bot}$ is an $n$-dimensional linear subspace spanned by 
$\{ {}^{\sharp} \phi^{1, 1}_1,  \ldots,  {}^{\sharp} \phi^{1, n}_1\}$.
Since $\dim  ( \ker \phi^1_1 )^{\bot} <\infty$,  the abstract Wiener space splits into two;
$ {\cal W} =    \overline{\ker \phi^1_1}^{ \| \,\cdot\,\|_{p -var} }  \oplus  ( \ker \phi^1_1 )^{\bot}$.
The projection $\pi$ naturally extends to the one from 
${\cal W}$ onto $\overline{\ker \phi^1_1}^{ \| \,\cdot\,\|_{p -var} } $,
which is again denoted by the same symbol.
There exist Gaussian measures $\mu_1$ and $\mu_2$ such that 
$(  \overline{\ker \phi^1_1}^{ \| \,\cdot\,\|_{p -var} }, \ker \phi^1_1, \mu_1 )$
and 
$( (\ker \phi^1_1)^{ \bot} ,  (\ker \phi^1_1)^{ \bot}  , \mu_2)$ are abstract Wiener spaces.
Naturally, $\mu_1 = \pi_* \mu$, $\mu_2 = \pi_*^{\bot} \mu$
and $\mu = \mu_1 \times \mu_2$ (the product measure).
One may think $\mu_1$ is the definition of the conditional measure 
${\mathbb P}[   \,\cdot\, |~ \phi^1_1  =0 ]$ in {\bf (A3)'} above.

Therefore, {\bf (A3)'} is equivalent to the following;
\begin{equation}\label{equia3.ineq}
{\mathbb E}[  \exp (  \la\bar\nu ,  \phi^2_1  ( \pi {\bf w}) \ra ) ]
< \infty.
\end{equation}
Precisely, $\pi {\bf w} 
:= {\cal L} ( \pi w)
=
\lim_{m \to \infty} {\cal L} \bigl( (\pi w)(m) \bigr)$.
Here,
${\cal L}$ stands for the rough path lift map.
Now we will see 
that $\pi {\bf w}$ is well-defined and has nice properties.


%
%
Note that $\phi_1^{1 } (k) = D_k \phi_1^{0 } (\bar\gamma)$
and recall (\ref{detuki.eq}), (\ref{Lmul2.eq}).
Hence, $\{  \phi^{1, 1}_1,  \ldots,  \phi^{1, n}_1\}$ 
are of rank $n$ in ${\cal H}^*$.
Let $C$
be the positive symmetric matrix in (\ref{detuki.eq}) 
and set  $K =(K_{ij}) = C^{-1}$,
$M =(M_{ij}) = C^{-1/2}$, which are again positive symmetric.
Then we have 
\begin{equation}
\pi w = w - \sum_{j=1}^n  \la w,  \sum_{l' =1}^n    M_{jl}   {}^{\sharp}\phi_1^{1, l}  \ra    
 \sum_{l' =1}^n M_{jl'}  {}^{\sharp}\phi_1^{1, l'} 
 =
 w - 
 \sum_{ l, l' =1}^n    K_{ll'}  \phi_1^{1, l} (w) \cdot  {}^{\sharp}\phi_1^{1, l'}.
 \label{proj1.eq}
\end{equation}  
This projection also works in $q$-variational setting.
Note that 
the second term on the right hand side is  ${\cal H}$-valued.
Therefore, the lift of $\pi w$
is actually a Young translation of ${\bf w}$ by 
$\sum_{ l, l'}  K_{l l'}  \phi_1^{1, l} (w) \cdot  {}^{\sharp}\phi_1^{1, l'}$.
It also holds that
$\pi {\bf w} = \lim_{m \to \infty} {\cal L} \bigl[ \pi ( w(m) )\bigr]$.
For $k, k' \in C^{q -var}_0 ([0,1], {\mathbb R}^d)$, we set
\begin{align}
 {\cal A} (k, k') 
&= 
\frac12 \hat{J}_1   \int_0^1 \hat{J}_t^{-1}
\{
\nabla \sigma ( \phi^0_t)\la \phi^1_t (k'),  dk_t \ra 
+
\nabla \sigma ( \phi^0_t)\la \phi^1_t (k),  dk'_t \ra 
\}
\nn\\
&\quad 
+ \frac12
 \hat{J}_1 \int_0^1 \hat{J}_t^{-1}  
  \nabla^2 \sigma ( \phi^0_t)\la \phi^1_t (k),   \phi^1_t (k'), d\bar\gamma_t \ra
 \end{align}
%
%
%
%
and $\hat{{\cal A}} (k,k')= \la \bar\nu,  {\cal A}  (\pi k, \pi k') \ra$.
Here, $\hat{J}  = \hat{J} (\bar\gamma)$ 
and $\phi^0 = \phi^0  (\bar\gamma)$
for brevity.
Then, $\hat{\cal A}$ is  a symmetric bounded bilinear mapping  form 
on  ${\cal H} \times {\cal H}$.
%
%
%
%
%
%
Notice that
\begin{equation}\label{Aphi2.eq}
{\cal A} (k,k) = \phi^2_1(k) -  \frac12 \delta^{H, 1/2} \cdot
\hat{J}_1 \int_0^1 \hat{J}_t^{-1}     b ( \phi^0_t)  dt,
\end{equation}
where $\delta^{H, 1/2} =1$ if $H =1/2$ and $\delta^{H, 1/2} =0$ if otherwise.
Therefore, 
$\hat{\cal A} (k,k) = \phi^2_1(k)  + {\rm (const)}$.

Now we will see that {\rm (i)}~ 
$\hat{\cal A}$ is actually Hilbert-Schmidt
and 
{\rm (ii)}~ $\phi^2_1 (\pi {\bf w}) \in {\cal C}_2 \oplus {\cal C}_0$
whose ${\cal C}_2$-component corresponds to $\hat{\cal A}$, 
that is,  $\phi^2_1 (\pi {\bf w}) = \Xi_{ \hat{\cal A} } (w)  + {\rm (const)}$.
Here,  
${\cal C}_j$ denotes the $j$th homogeneous Wiener chaos of order $j$
and 
$\Xi_{ {\cal B} }$  denotes the element in  $ {\cal C}_2$ 
which unitarily corresponds to a symmetric Hilbert-Schmidt  
bilinear form ${\cal B}$.


For $m \in {\mathbb N}$,  
set $\hat{{\cal A}}_m (k,k')= \la \bar\nu,  {\cal A}  ( (\pi k)(m), (\pi k')(m)) \ra$.
The corresponding bounded self-adjoint operator on  ${\cal H}$  is denoted by $\hat{A}_m$.
This bilinear form extends to a bounded bilinear form on 
${\cal W} \times {\cal W}$.
Hence, by Goodman's theorem (see Theorem 4.6, p. 83, \cite{kuo}),
it is of trace class (and consequently Hilbert-Schmidt).
$\hat{{\cal A}}_m (w,w) =  \Xi_{\hat{{\cal A}}_m } (w)  + {\rm Trace} (\hat{A}_m)$.
As a result, 
$\phi^2_1 ( (\pi  w)(m)  ) = \Xi_{\hat{{\cal A}}_m } (w)  +  s_m$, 
where the constant $s_m$ may depend on $m$.


%
By a straight forward rough path calculation as in Section 5, Inahama \cite{inaaop}, 
we can prove that $\phi^2_1 ((\pi w)(m) )$ converges to 
$\phi^2_1 (\pi {\bf w})$ in $L^2 (\mu)$.
(In Inahama \cite{inaaop},  the convergence 
$\phi^2_1 ( w(m) ) \to \phi^2_1 ( {\bf w} ) $ as $m \to \infty$ is shown.
We can  modify that proof, 
since the effect of the projection $\pi$ appears as Young translation 
as we have already seen.)
Hence,  both $\Xi_{\hat{{\cal A}}_m }$ and $s_m$ 
converge in  $ {\cal C}_2$ and $ {\cal C}_0$, respectively.
By the unitary correspondence, 
there exists  a symmetric Hilbert-Schmidt  bilinear form ${\cal B}$
such that $\hat{{\cal A}}_m \to {\cal B}$ as $m \to \infty$
in Hilbert-Schmidt norm.
From a basic property of Young integral, we see that 
$\hat{{\cal A}}_m (k,k') \to \hat{{\cal A}} (k,k') $ as $m \to \infty$ 
for each fixed $k, k' \in {\cal H}$.
Thus we have shown  {\rm (i)} and {\rm (ii)} above.


Exponentially integrability of  quadratic Wiener functionals 
is well-known.
(\ref{equia3.ineq}) is equivalent to 
${\mathbb E}[  \exp (  \Xi_{\hat{\cal A}  }  ) ]< \infty$, which in turn is equivalent to 
$\sup {\rm Spec }(\hat{A}) <1/2$.
Since the inequality is strict, 
there exists $\rho >1$ such that $\sup {\rm Spec }(\rho \hat{A}) <1/2$,
which is equivalent to ${\mathbb E}[  \exp ( \rho \hat\Xi_{\hat{\cal A}  }   ) ]< \infty$.
Summing it up, we have seen that  {\bf (A3)'} is
equivalent to the following;
\begin{equation}\label{equia7.ineq}
{\mathbb E}[  \exp ( \rho \la\bar\nu ,  \phi^2_1( \pi {\bf w} ) \ra ) ]
< \infty
\qquad \mbox{  for some $\rho >1$.}
\end{equation}
%
%



Let us check here that {\bf (A3)} and {\bf (A3)'}
are equivalent under {\bf (A1)}, {\bf (A2)}.
\begin{pr}\label{a3equi.pr}
Under {\bf (A1)} and {\bf (A2)}, the two conditions {\bf (A3)} and {\bf (A3)'}
are equivalent.
\end{pr}

\Proof
As is explained above, ${\bf (A3)'}$ is equivalent to $\sup {\rm Spec }(\hat{A}) <1/2$.
Keep in mind that the only accumulation point of ${\rm Spec }(\hat{A})$ is $0$,
since $\hat{A}$ is Hilbert-Schmidt.
Let $(- \ve_0, \ve_0) \ni u \mapsto f(u) \in K_a^{a'}$ be a smooth curve in $K_a^{a'}$
such that $f(0) = \bar\gamma$ and $f^{\prime} (0)  \neq 0$ as in {\bf (A3)}.
Then, a straight forward calculation shows that
\begin{align}
\lefteqn{
\frac{d^2}{du^2} \Big|_{u=0} \frac{ \|   f(u) \|^2_{{\cal H}} }{2}
=
\frac{d^2}{du^2} \Big|_{u=0}
\Bigl(  \frac{ \|   f(u) \|^2_{{\cal H}} }{2}  - \la \bar\nu, \phi^0_1(f_u) -a' \ra  \Bigr)
}
\nn\\
&
=
\|   f^{\prime}(0) \|^2_{{\cal H}}+ \la  f^{\prime\prime}(0),  \bar\gamma \ra_{{\cal H}}
 - \bigl\la \bar\nu, D\phi^0_1(\bar\gamma) \la f^{\prime\prime}(0) \ra \bigr\ra 
 - \bigl\la \bar\nu, D^2\phi^0_1(\bar\gamma) \la f^{\prime}(0), f^{\prime}(0) \ra \bigr\ra 
    \nn\\
&
=
\|   f^{\prime}(0) \|^2_{{\cal H}}
-  \bigl\la \bar\nu, D^2\phi^0_1(\bar\gamma) \la \pi f^{\prime}(0), \pi f^{\prime}(0) \ra \bigr\ra
\nn\\
&
=
\|   f^{\prime}(0) \|^2_{{\cal H}}
-  2 \bigl\la \bar\nu, \psi \la \pi f^{\prime}(0), \pi f^{\prime}(0) \ra \bigr\ra,
 \label{2kaibi.eq}
    \end{align}
where we used (\ref{Lmul1.eq})--(\ref{Lmul2.eq}) and the fact that $f^{\prime}(0)$
is tangent to the submanifold $K_a^{a'}$.
Since $f^{\prime}(0)$ can be any non-zero vector $h$ such that $\pi h =h$, 
we see from (\ref{2kaibi.eq}) that {\bf (A3)} is equivalent to 
\[
\bigl\la \bar\nu, \psi \la \pi f^{\prime}(0), \pi f^{\prime}(0) \ra \bigr\ra
< 
  \frac12  \| h \|^2_{{\cal H}}
\qquad
\qquad
(h \in {\cal H} \setminus \{0\}),
\]
which in turn is equivalent to
$\sup {\rm Spec } (\hat{A}) <1/2$.
\QED


The following is a key technical lemma. 
Roughly speaking, it states that  restricted on a sufficiently small subset,
 $\exp(\la \bar\nu , r_{\ve, 1}^{2} \ra /\ve^2) \in \cup_{1< q <\infty} L^q$ uniformly in $\ve$.
\begin{lm}\label{lm.integ3}
Assume {\bf (A1)}, {\bf (A2)} and {\bf (A3)}.
Then, there exists $\rho_1 >1$ and $\eta >0$ such that
\[
\sup_{0 < \ve \le 1}
{\mathbb E} \bigl[
\exp \bigl(
   \rho_1  \la \bar\nu ,  r_{\ve,1}^2   \ra /\ve^2
           \bigr) 
I_{ U_{\eta} }  (\ve {\bf w})
%
%
I_{  \{   | r^1_{\ve, 1}   /\ve | \le \eta_1  \}}
\bigr]  <\infty
\]
for any $\eta_1 >0$.
\end{lm}


\Proof
By Lemma \ref{lm.integ1} and the relation $r_{1,\ve}^{2}  /\ve^2 = \phi^2_1 +r_{\ve,1}^{2+}  /\ve^2 $,
it is sufficient to show that
\begin{equation}
\sup_{0 < \ve \le 1}
{\mathbb E} \bigl[
\exp \bigl(
    \rho_1  \la \bar\nu ,   \phi^2_1 \ra 
           \bigr) 
%
%
I_{ U_{\eta} }  (\ve {\bf w})
I_{  \{   | r_{\ve,1}^{1}   /\ve | \le \eta_1  \}}
\bigr]  <\infty.
\label{integ3suff.eq}
\end{equation}

%
%
%
%
%
%

Then,  from (\ref{proj1.eq}) and (\ref{Aphi2.eq}) we have
\begin{align}
 \phi^2_1 ({\bf w})  
 &=  
\lim_{m\to \infty}  \phi^2_1 (w(m))   
=
\lim_{m\to \infty}  {\cal A} (w(m), w(m)) -{\rm (const)}
 \nn\\
    & =\phi^2_1 ( \pi {\bf w})  + 2    \sum_{j,j'}   \phi_1^{1,j} (w)
K_{j j'} 
\cdot   {\cal A}  \la w,   {}^{\sharp} \phi_1^{1,j'} \ra
   \nn\\
      & \qquad  + 
          \sum_{j,j' , k, k'}  \phi_1^{1,j} (w) \phi_1^{1,k} (w) K_{j j'} K_{k k'}
           \cdot   {\cal A}  \la   {}^{\sharp} \phi_1^{1,j'}  ,   {}^{\sharp} \phi_1^{1,k'}  \ra
=: Z_1+ Z_2+ Z_3.
                                     \label{proj2.eq}
                                        \end{align}
Note that 
${\cal A}  \la w,   {}^{\sharp} \phi_1^{1,j'} \ra$ 
and 
$ {\cal A}  \la   {}^{\sharp} \phi_1^{1,j'}  ,   {}^{\sharp} \phi_1^{1,k'}  \ra$
 are well-defined 
as Young integrals.

Exponential integrability of the first term $Z_1$ 
on the right hand side of (\ref{proj2.eq}) is given in (\ref{equia7.ineq}).
So, we estimate the second term $Z_2$.
Since $\ve \phi_1^{1} (w) = r_{\ve, 1}^{2} ({\bf w})-  r_{\ve, 1}^{1} ({\bf w})$
and $|{\cal A}  \la w,   {}^{\sharp} \phi_1^{1,j'} \ra|
\lesssim \|w\|_{p -var}$, 
we have
\begin{align}
|\phi_1^{1,j} (w) {\cal A}  \la w,   {}^{\sharp} \phi_1^{1,j'} \ra|
&\le
c_1 \Bigl\{
\Bigl|\frac{  r_{\ve, 1}^{2}  ({\bf w})  }{\ve}\Bigr| 
+ \Bigl| \frac{ r_{\ve, 1}^{1} ({\bf w}) }{\ve}  \Bigr|   \Bigr\}
\|w\|_{p -var}  
\nn\\
&\le
c_1 \Bigl\{  \Bigl|\frac{ c'r_{\ve, 1}^{2} ({\bf w})  }{\ve}\Bigr|^2 
+ \frac{\|w\|_{p -var}^2 }{4 c'^2}   \Bigr\}
+
c_1 \Bigl| \frac{  r_{\ve, 1}^{1} ({\bf w}) }{\ve}  \Bigr|    \|w\|_{p -var}  
\nn
\end{align}
for any $c' >0$.

Set $c_2 =2c_1n^2 \sup_{j,j'} |K_{j,j'}|$ 
and let $M>0$. 
Then, by H\"older's inequality,
\begin{align}
{\mathbb E} \bigl[
e^{M |Z_2| }
%
%
I_{ U_{\eta} }  (\ve {\bf w})
I_{  \{   | r_{\ve, 1}^{1}   /\ve | \le \eta_1  \}}
\bigr]  
\le
{\mathbb E} \bigl[
\exp \bigl( 3M c_2 c'^2  |r_{\ve,1}^{2}   /\ve|^2 \bigr)
%
%
I_{ U_{\eta} }  (\ve {\bf w})
\bigr]^{1/3} 
\nn\\
 \quad  \times
{\mathbb E} \bigl[
e^{  3M c_2 \|w\|_{p -var}^2 /(4 c')}
\bigr]^{1/3} 
{\mathbb E} \bigl[
e^{  3M c_2 \eta_1 \| w \|_{p -var}  }
\bigr]^{1/3}.
\nn
\end{align}
For any $M>0$ and $\eta_1 >0$, the third factor is integrable.
If $c'$ is chosen sufficiently large, 
then the second factor is also integrable by Fernique's theorem.
By Lemma \ref{lm.integ2},  there exists $\eta >0$ such that $\sup_{\ve}$ of the first 
factor is finite and, hence, 
\begin{equation}\label{proj3.eq}
\sup_{0< \ve \le 1}
{\mathbb E} \bigl[
e^{M |Z_2| }
%
%
I_{ U_{\eta} }  (\ve {\bf w})
I_{  \{   | r^1_{ \ve, 1}   /\ve | \le \eta_1  \}}
\bigr]  
< \infty.
\end{equation}

Since $\phi_1^{1,j} (w) \phi_1^{1,k} (w) 
= \ve^{-1} \{ r^2_{ \ve,1} ({\bf w})^j- r_{\ve, 1}^{1} ({\bf w})^j \} \phi_1^{1,k} (w)$,
we can deal with $Z_3$ in the same way.
For any $M>0$ and $\eta_1 >0$,  there exists $\eta >0$ such that
\begin{equation}\label{proj4.eq}
\sup_{0< \ve \le 1}
{\mathbb E} \bigl[
e^{M |Z_3| }
%
%
I_{ U_{\eta} }  (\ve {\bf w})
I_{  \{   | r_{1, \ve}^{1}   /\ve | \le \eta_1  \}}
\bigr]  
< \infty.
\end{equation}

Let $\rho >1$ be as in (\ref{equia7.ineq}).
Set $\rho_1=(1+ \rho)/2 >1$, $s =2\rho/ (1+ \rho) >1$, and $1/s +1/s' =1$.
Then, from H\"older's inequality and (\ref{equia7.ineq}), (\ref{proj2.eq})--(\ref{proj4.eq}),
we can easily see that
\begin{align}
\lefteqn{
{\mathbb E} \bigl[
\exp \bigl(
    \rho_1  \la \bar\nu ,   \phi^2_1 \ra 
           \bigr) 
%
%
I_{ U_{\eta} }  (\ve {\bf w})
I_{  \{   | r^1_{1, \ve}   /\ve | \le \eta_1  \}}
\bigr]  
}
\nn\\
&\le
{\mathbb E} \bigl[  \exp \bigl(
    \rho \la \bar\nu ,   \phi^2_1 \circ \pi \ra 
           \bigr) \bigr]^{1/s}
\prod_{i=1}^2
{\mathbb E} \bigl[
e^{ 2q' \rho_1 |\bar\nu| |Z_i| }
I_{ U_{\eta} }  (\ve {\bf w})
I_{  \{   | r^1_{1, \ve}   /\ve | \le \eta_1  \}}
%
%
\bigr]^{1/(2s')}.
\nn
\end{align}
From this, (\ref{integ3suff.eq}) is immediate.
This completes the proof.
\QED

\subsection{Proof of off-diagonal short time asymptotics}

In this subsection we prove Theorem \ref{thm.MAIN.off}, namely,
off-diagonal short time asymptotics
of the density of the solution $(y_t) = (y_t (a))$ 
of RDE (\ref{main.ygODE.eq})
driven by fractional Brownian rough path
 ${\bf w}$ with $1/3 <H \le 1/2$ under Assumptions {\bf (A1)}--{\bf (A3)}. 

%
%

First, let us calculate the kernel $p(t, a, a')$.
Take  $\eta >0$ as in Lemma \ref{lm.integ3}.
Then, we see
\begin{align}
p(\ve^{1/H}, a, a') 
&=
{\mathbb E} \bigl[  
\delta_{a'} ( y_1^{\ve})   
 \bigr]
 \nn\\
 &= 
   {\mathbb E} \bigl[  
\delta_{a'} ( y_1^{\ve} )   \chi_{\eta} (\ve, w)
 \bigr]
   + 
    {\mathbb E} \bigl[  
\delta_{a'} ( y_1^{\ve} )  \bigl\{ 1- \chi_{\eta} (\ve, w)  \bigr\}
 \bigr]
 = : I_1 +I_2.
 \nn
  \end{align}
As we have shown in Lemma \ref{lm.ldpcut}, the second term $I_2$ on the right hand side 
does not contribute to the asymptotic expansion.
So, we have only to calculate the first term $I_1$.
By Cameron-Martin formula, 
\[
I_1 
=
   {\mathbb E} \bigl[  
   \exp \bigl(  -\frac{ \|  \bar\gamma \|^2_{{\cal H}}}{2\ve^2}   - \frac{1}{\ve} \la \bar\gamma, w \ra \bigr)
\delta_{a'} ( \tilde{y}_1^{\ve} )   \chi_{\eta} (\ve, w + \frac{\bar\gamma}{\ve})
 \bigr].
 \]
Recall that $\la \bar\gamma, w \ra = \la \bar\nu , \phi^1_1 (w) \ra$ for all $w$.
Hence,  we have
\begin{align}
I_1  
&=
 \exp \bigl(  -\frac{ \|  \bar\gamma \|^2_{{\cal H}}}{2\ve^2}  \bigr)
  {\mathbb E} \bigl[  
   \exp \bigl(
    - \frac{1}{\ve} \la \bar\nu , \phi^1_1  \ra
           \bigr)
\delta_{a'} ( a' +\ve \phi_1^1    +  r_{\ve, 1}^{2} )   \chi_{\eta} (\ve, w + \frac{\bar\gamma}{\ve})
 \bigr]
    \nn\\
      &=  \frac{1}{\ve^n}  \exp \bigl(  -\frac{ \|  \bar\gamma \|^2_{{\cal H}}}{2\ve^2}  \bigr)
  {\mathbb E} \bigl[  
   \exp \bigl(
    - \frac{1}{\ve} \la \bar\nu , \phi^1_1  \ra
           \bigr)
\delta_{0} (  \phi_1^1 + \ve^{-1}   r_{\ve, 1}^{2} )   
\chi_{\eta} (\ve, w + \frac{\bar\gamma}{\ve})
 \bigr]
    \nn\\
      &=    
        \frac{1}{\ve^n}  \exp \Bigl(  -\frac{ \|  \bar\gamma \|^2_{{\cal H}}}{2\ve^2}  
       \Bigr)
    %
    %
  {\mathbb E} \bigl[  
   \exp \bigl(
      \la \bar\nu ,  r_{\ve, 1}^{2}   \ra /\ve^2
           \bigr)
\delta_{0} (  \phi_1^1 + \ve^{-1} r_{\ve, 1}^{2} )   \chi_{\eta} (\ve, w + \frac{\bar\gamma}{\ve})
 \bigr]
    \nn\\
      &= 
     \frac{1}{\ve^n}  \exp \Bigl(  -\frac{ \|  \bar\gamma \|^2_{{\cal H}}}{2\ve^2}  
      \Bigr)
        {\mathbb E} \bigl[      F(\ve, w)
\delta_{0} \bigl(  \frac{ \tilde{y}_1^{\ve} - a'}{\ve } \bigr) 
 \bigr],
 \nn
       \end{align}
where 
\begin{align}
F(\ve, w) =  \exp \bigl(  \ve^{-2}
      \la \bar\nu ,  r^2_{ \ve,1}   \ra 
           \bigr) 
           \chi_{\eta} (\ve, w + \frac{\bar\gamma}{\ve})       
\psi \Bigl(  \frac{1}{\eta_1^2}  \Bigl| \frac{\tilde{y}_1^{\ve} - a'}{\ve}   \Bigr|^2 \Bigr)
\label{Fve.def}
\end{align}
for any positive constant $\eta_1$. 
Here, $\psi$ is the cut-off function introduced in Subsection \ref{subsec.ldp}.
It is easy to see that {\rm (i)}~
$\chi_{\eta} (\ve, w +  \bar\gamma /\ve)  $ and its derivatives vanish 
outside $\{  w~|~ \ve {\bf w} \in U_{\eta} \}$
and  {\rm (ii)}~
$\psi \bigl( \eta_1^{-2}  \bigl| (\tilde{y}_1^{\ve} - a' )/\ve   \bigr|^2  \bigr)$ 
and its derivatives vanish 
outside $\{  |r^1_{\ve,1 } /\ve | \le \eta_1 \}$.
Hence, by Lemma \ref{lm.integ3}, $F(\ve, w) \in \tilde{\bf D}_{\infty}$ and 
$F(\ve, w) =O(1)$ with respect to that topology.
Since $\delta_{0} (  ( \tilde{y}_1^{\ve} - a')/\ve  )$ admits an asymptotic expansion
in $\tilde{\bf D}_{-\infty}$,
the problem reduces to whether 
$F(\ve, w) $  admits an asymptotic expansion
in $\tilde{\bf D}_{\infty}$.


\begin{lm}\label{lm.sikihen1}
Assume {\bf (A1)}--{\bf (A3)}.
For any $M \in {\mathbb N}$,  we have
\[
  {\mathbb E} \bigl[      F(\ve, w)
\delta_{0} \bigl(  \frac{ \tilde{y}_1^{\ve} - a'}{\ve}  \bigr) 
 \bigr]
=
  {\mathbb E} \bigl[      F(\ve, w)  \psi ( | \phi_1^1 / \eta_1|^2)
\delta_{0} \bigl(  \frac{ \tilde{y}_1^{\ve} - a'}{\ve}  \bigr) 
 \bigr]
+
O(\ve^M)
\]
as $\ve \searrow 0$.
\end{lm}

\Proof
By using Taylor expansion for $\psi$, we see that, for given $M$,
there exist $m \in {\mathbb N}$ and $G_j (\ve, w) \in {\bf D}_{\infty} ~(1 \le j \le m)$
such that
\begin{align}
\psi \Bigl(  \frac{1}{\eta_1^2}  \Bigl| \frac{\tilde{y}_1^{\ve} - a'}{\ve}   \Bigr|^2 \Bigr)
&=
\psi \Bigl( \bigl| \frac{\phi_1^1}{\eta_1}  \bigr|^2 \Bigr)
+ \sum_{j=1}^m \psi^{(j)} \Bigl( \bigl| \frac{\phi_1^1}{\eta_1} \bigr|^2   \Bigr) G_j (\ve, w)
+O(\ve^M)
\label{siki2.eq}
\end{align}
in ${\bf D}_{\infty}$ as $\ve \searrow 0$.
$G_j (\ve, w) =O(1)$,  but its explicit form is not important.
Note that $\psi^{(j)} (| \phi_1^1/\eta_1|^2 ) T(\phi_1^1)=0$ 
if $j \ge 1$ and ${\rm supp} (T) \subset \{ a \in {\mathbb R}^n~|~|a| <  \eta_1 /2\}$.

By Proposition \ref{pr.nondeg} and 
Watanabe's asymptotic theory in \cite{wa, iwbk},
$\delta_{0} (  ( \tilde{y}_1^{\ve} - a')/\ve )$ admits an asymptotic expansion
in $\tilde{\bf D}_{-\infty}$ as follows.
As before, we set $\{ 0= \nu_0 < \nu_1 < \nu_2 <\cdots\}$ to be all the elements of $\Lambda_3$
in increasing order. 
For given $M$, let $l \in {\mathbb N}$ be the smallest integer such that $M \le \nu_{l+1}$.
Then, for some 
$\Phi_{\nu_j}  \in \tilde{\bf D}_{-\infty} ~(1 \le j \le l)$,
it holds that
\begin{align}
\delta_{0} (  ( \tilde{y}_1^{\ve} - a')/\ve )
=
\delta_{0} ( \phi_1^1 )
+
\ve^{\nu_1} \Phi_{\nu_1} +\cdots + \ve^{\nu_l} \Phi_{\nu_l}
+O(\ve^{ \nu_{l+1}})
\label{siki3.eq}
\end{align}
in $\tilde{\bf D}_{-\infty}$ as $\ve \searrow 0$.
Here, $\Phi_{\nu_j}$ is a finite linear combination of terms of the form
\[
\partial^{\beta} \delta_0 (\phi_1^1) \times 
\{
\mbox{a polynomial of the components of $\phi_1^{\kappa_i}$\,'s}
\},
\]
where $\beta$ stands for a multi-index.
Hence, $\psi^{(j')} (| \phi_1^1/\eta_1|^2 )  \Phi_{\nu_j}$ vanish for all $j, j'$.

Now, using (\ref{siki2.eq}) and (\ref{siki3.eq}), we prove the lemma. 
\begin{align}
\lefteqn{ 
 {\mathbb E} \bigl[      F(\ve, w)
\delta_{0} (  ( \tilde{y}_1^{\ve} - a')/\ve ) 
 \bigr]
}
\nn\\
&=
 {\mathbb E} \bigl[      F(\ve, w)
\psi \Bigl(  \frac{1}{\eta_1^2}  \Bigl| \frac{\tilde{y}_1^{\ve} - a'}{\ve}   \Bigr|^2 \Bigr)
\delta_{0} (  ( \tilde{y}_1^{\ve} - a')/\ve ) 
 \bigr]
\nn\\
&=
 {\mathbb E} \bigl[      F(\ve, w)
\psi ( | \phi^1_1/\eta_1|^2 )
\delta_{0} (  ( \tilde{y}_1^{\ve} - a')/\ve ) 
 \bigr]
\nn\\
&  \quad
+  {\mathbb E} \bigl[      F(\ve, w)
\Bigl(  
\sum_{j=1}^m \psi^{(j)} \Bigl( \bigl| \frac{\phi_1^1}{\eta_1} \bigr|^2   \Bigr) G_j (\ve, w)
 \Bigr)
\delta_{0} (  ( \tilde{y}_1^{\ve} - a')/\ve ) 
 \bigr]
+
O(\ve^M)
\nn\\
&=
 {\mathbb E} \bigl[      F(\ve, w)
\psi ( | \phi^1_1/\eta_1|^2 )
\delta_{0} (  ( \tilde{y}_1^{\ve} - a')/\ve ) 
 \bigr]
\nn\\
&  \quad
+  {\mathbb E} \bigl[      F(\ve, w)
\Bigl(  
\sum_{j=1}^m \psi^{(j)} \Bigl( \bigl| \frac{\phi_1^1}{\eta_1} \bigr|^2   \Bigr) G_j (\ve, w)
 \Bigr)
\bigl(  \delta_{0} ( \phi_1^1 )
+
\cdots + \ve^{\nu_l} \Phi_{\nu_l} \bigr)
 \bigr]
+
O(\ve^M)
\nn\\
&=
  {\mathbb E} \bigl[      F(\ve, w)  \psi ( | \phi_1^1 / \eta_1|^2)
\delta_{0} (  ( \tilde{y}_1^{\ve} - a')/\ve ) 
 \bigr]
+
O(\ve^M).
\nn
\end{align}
Thus, we have shown the lemma.
\QED


Set 
$
\Lambda'_2 = \{ \kappa -2 ~|~ \kappa \in \Lambda_1 \setminus \{0, 1\} \}$.
If $H \neq 1/2$,
then 
$\Lambda'_2 =
\{ 0 <  H^{-1} -2 < 1  <\cdots \}.
$
Next we set 
$
\Lambda'_3 = 
\{  a_1 +a_2 + \cdots + a_m ~|~  \mbox{$m \in {\bf N}_+$ and $a_1 ,\ldots, a_m \in \Lambda'_2$} \}$.
In the following lemma, 
$\{ 0=\xi_0 <\xi_1<\xi_2 <\cdots \}$ stands for all the elements of 
$\Lambda'_3$ in increasing order.

Note that the following lemma
does not claim $F_{k+1} (\ve, w)  = O(\ve^{\xi_{k+1}})$, 
but it claims 
$F_{k+1} (\ve, w) T ( \phi^1_1 ) = O(\ve^{\xi_{k+1}})$
if  $T \in {\cal S}^{\prime} ({\mathbb R}^n)$
is for example of the form $\partial^{\beta} \delta_0$.

%
%
%

\begin{lm}\label{lm.asyFve}
Assume {\bf (A1)}--{\bf (A3)} and let $F(\ve, w) \in \tilde{\bf D}_{\infty}$ as in (\ref{Fve.def}).
Then, for every $k=1,2,3, \ldots$,
\begin{align}
\lefteqn{
F(\ve, w) \psi (  |  \phi^1_1 (w)  / \eta_1 |^2 )
}
\nn\\
 &= \exp \bigl(
      \la \bar\nu ,  \phi^2_1 ({\bf w})
           \bigr) 
\psi (  |  \phi^1_1 (w)  / \eta_1 |^2 )^2
 \{1+ \ve^{\xi_1} K_{\xi_1}(w)  +\cdots + \ve^{\xi_k} K_{\xi_k} (w)\}
+
F_{k+1} (\ve, w),
\nn
\end{align}
where $F_{k+1} (\ve, w) \in \tilde{\bf D}_{\infty}$ satisfies that
\[
F_{k+1} (\ve, w) T ( \phi^1_1 ) = O(\ve^{\xi_{k+1}})  \qquad
\mbox{ in ${\bf D}_{-\infty}$ as $\ve \searrow 0$}
\]
for any $T \in {\cal S}^{\prime} ({\mathbb R}^n)$ with ${\rm supp}(T) 
\subset \{ a\in {\mathbb R}^n ~|~ |a| \le \eta_1/2 \}$.
Moreover, $K_{\xi_j} \in {\bf D}_{\infty} ~(j=1,2,\ldots)$ are determined by the following 
formal expansion ($\kappa_3 =H^{-1}$ if $H \neq 1/2$);
\begin{align}
\sum_{m=0}^{\infty}  \frac{
\la  \bar\nu,  r_{\ve, 1}^{\kappa_3 }  /\ve^2\ra^m }{m!}
&=
\sum_{m=0}^{\infty}
\frac{1}{m!}
\Bigl\{
  \ve^{\kappa_3 -2} \la  \bar\nu,  \phi_1^{\kappa_3} \ra
+  \ve^{\kappa_4 -2} \la  \bar\nu, \phi_1^{\kappa_3}  \ra +\cdots   
\Bigr\}^m
\nn\\
&=
 1+ \ve^{\xi_1} K_{\xi_1} + \ve^{\xi_2}  K_{\xi_2}  +\cdots.
\nn
\end{align}
\end{lm}

%
%
%

\Proof
Let $\rho_1 >1$ be as in Lemma \ref{lm.integ3}.
First we show that, for any $\eta_1 >0$, 
\begin{equation}
{\mathbb E} \bigl[
\exp \bigl(
    \rho_1  \la \bar\nu ,   \phi^2_1 \ra 
           \bigr) 
I_{  \{   | \phi_1^1 | \le \eta_1  \}}
\bigr]  <\infty.
\label{siki1.eq}
\end{equation}
We can choose a subsequence $\{ \ve_k \}$
such that, as $k \to \infty$, $\ve_k \searrow 0$ and 
$R_1^{1, \ve_k} / \ve_k  \to \phi^1_1$ a.s.
To prove (\ref{siki1.eq}), we apply Fatou's lemma
to  (\ref{integ3suff.eq}) with $\eta_1$ replaced by $2\eta_1$.
\begin{align}
\infty &> \liminf_{k \to \infty}
{\mathbb E} \bigl[
\exp \bigl(
    \rho_1  \la \bar\nu ,   \phi^2_1 \ra 
           \bigr) 
%
%
I_{ U_{\eta} }  ( \ve_k {\bf w})
I_{  \{   | r_{\ve_k, 1}^{1}   /\ve_k | \le 2\eta_1  \}}
\bigr] 
\nn\\
&\ge 
{\mathbb E} \bigl[
\exp \bigl(
    \rho_1  \la \bar\nu ,   \phi^2_1 \ra 
           \bigr) 
\liminf_{k \to \infty}
I_{  \{   | r^1_{ \ve_k, 1}   /\ve_k | \le 2\eta_1  \}}
\bigr] 
\ge 
{\mathbb E} \bigl[
\exp \bigl(
    \rho_1  \la \bar\nu ,   \phi^2_1 \ra 
           \bigr) 
I_{  \{   |\phi^1_1 | \le \eta_1  \}}
\bigr]. 
\nn
\end{align}
From (\ref{siki1.eq}), it is easy to check that 
$\exp \bigl(
      \la \bar\nu ,  \phi^2_1  \ra
           \bigr) 
\psi (  |  \phi^1_1  / \eta_1 |^2 ) \in \tilde{\bf D}_{\infty}$.

Now we expand $\exp ( \la \bar\nu ,  r_{\ve, 1}^{2}   \ra /\ve^2 ) 
= \exp (    \la \bar\nu ,  \phi^2_1  \ra ) \exp ( \la \bar\nu ,  r_{\ve,1}^{\kappa_3 }   \ra /\ve^2 ) $ in $\ve$.
Set $Q_{l+1} : {\mathbb R} \to {\mathbb R}$ by
\[
Q_{l+1} (u) = e^u -\Bigl(  1+ u +\frac{u^2}{2!}  +\cdots +  \frac{u^l}{l!}  \Bigr) = 
u^{l+1} 
\int_0^1 \frac{(1-\theta)^l }{l!} e^{\theta u}  d\theta 
\qquad
(u \in {\mathbb R}).
\]
We will prove that, for sufficiently large $l \in {\mathbb N}$, as $\ve \searrow 0$,
\begin{equation}\label{siki4.eq}
e^{ \la \bar\nu ,  \phi^2_1  \ra }
Q_{l+1}( \la  \bar\nu, r_{\ve, 1}^{\kappa_3 } \ra / \ve^2 )   \chi_{\eta} (\ve, w + \frac{\bar\gamma}{\ve})  \psi (  |  \phi^1_1  / \eta_1 |^2 ) 
=
O( \ve^{\xi_{k+1}})
\quad
\mbox{in $\tilde{\bf D}_{\infty}$. }
\end{equation}
Note that $  \chi_{\eta} (\ve, w + \frac{\bar\gamma}{\ve})  =O(1)$ in ${\bf D}_{\infty}$ as $\ve \searrow 0$ by (\ref{chi_asy.eq}).
By Proposition \ref{pr.asyDinf},
$r_{\ve,1}^{\kappa_3 } / \ve^2 =O(\ve^{  \kappa_3 -2})$ in ${\bf D}_{\infty}$.
So, if $l +1 \ge \xi_{k+1}/  (\kappa_3 -2 )$, 
then 
$( \la \bar\nu, r_{\ve, 1}^{\kappa_3} \ra / \ve^2 )^{l+1} =O(\ve^{ \xi_{k+1}})$ in ${\bf D}_{\infty}$.
Therefore, in order to verify (\ref{siki4.eq}), it is sufficient to show that, as $\ve \searrow 0$,
\begin{equation}\label{siki5.eq}
\int_0^1 (1-\theta)^l  
e^{ \la     \bar\nu ,  \phi^2_1 +   \theta   r_{\ve, 1}^{\kappa_3} / \ve^2 \ra }  d\theta \cdot
 \chi_{\eta} (\ve, w + \frac{\bar\gamma}{\ve})  \psi (  |  \phi^1_1  / \eta_1 |^2 ) 
=
O( 1)
\quad
\mbox{in $\tilde{\bf D}_{\infty}$. }
\end{equation}
To verify the integrability of this Wiener functional, 
note that $e^{\theta u}  \le 1+ e^u$ for all $u \in {\mathbb R}$
and $0 \le \theta \le 1$.
This implies that the first factor on the left hand side of (\ref{siki5.eq})
 is dominated by $e^{ \la \bar\nu ,  \phi^2_1  \ra } + e^{\la \bar\nu, r_{\ve,1}^{2} \ra / \ve^2}$.
From Lemma \ref{lm.integ3}
and (\ref{siki1.eq}),  we see that the left hand side of  (\ref{siki5.eq}) is 
 $O(1)$ in any $L^r ~(1<r<\infty)$.
In the same way, 
the Malliavin derivatives of the left hand side of (\ref{siki5.eq}) are $O(1)$ in any $L^r$.

It is easy to see that, as $\ve \searrow 0$,
\begin{align}
\sum_{k=0}^l
\frac{  \{  \la  \bar\nu, r_{\ve, 1}^{\kappa_3} \ra  / \ve^2 \}^k  }{k!}  
=
 1+ \ve^{\xi_1} K_{\xi_1} +\cdots + \ve^{\xi_k} K_{\xi_k}  
 +
O( \ve^{\xi_{k+1}})
\quad
\mbox{in ${\bf D}_{\infty}$. }
\end{align}
From this and (\ref{chi_asy.eq}), we see that 
\begin{align}
\lefteqn{
F(\ve, w) \psi (  |  \phi^1_1 (w)  / \eta_1 |^2 )
}
\nn\\
 &= \exp \bigl(
      \la \bar\nu ,  \phi^2_1 ({\bf w})
           \bigr) 
\psi (  |  \phi^1_1 (w)  / \eta_1 |^2 )  
 \psi \Bigl(  \frac{1}{\eta_1^2}  \Bigl| \frac{\tilde{y}_1^{\ve} - a'}{\ve}   \Bigr|^2 \Bigr) \{1+ \ve^{\xi_1} K_{\xi_1}(w)  +\cdots + \ve^{\xi_k} K_{\xi_k} (w)\}
\nn\\
& \qquad\qquad 
 +
O( \ve^{\xi_{k+1}})
\quad
\mbox{in $\tilde{\bf D}_{\infty}$. }
\nn
\end{align}
Using (\ref{siki2.eq}), we finish the proof.
\QED


{\it Proof of the main theorem (Theorem \ref{thm.MAIN.off})}~
Now we prove our main theorem in this paper.
We set 
$$\Lambda_4 =\Lambda_3  +\Lambda'_3 =\{ \nu + \xi ~|~ \nu   \in  \Lambda_3  , \xi \in  \Lambda'_3\}.$$
We denote by $\{ 0 = \lambda_0 < \lambda_1 < \lambda_2 < \cdots \}$ all the elements of $\Lambda_4$  in increasing order.
It is no mystery why this index set  appears in the short time expansion of the kernel 
because, very formally speaking,  the problem reduces to finding asymptotic behavior  of
$ {\mathbb E}[  \exp(  \la \bar\nu, r_{\ve,1}^{2 } \ra  / \ve^2 )  
\cdot \delta_0 (r^1_{ \ve,1} / \ve) ]$, as we have seen.
Now, by (\ref{Fve.def}), 
Lemma \ref{lm.sikihen1}
Lemma \ref{lm.asyFve}, and (\ref{siki3.eq}), we can easily prove 
the asymptotic expansion in Theorem \ref{thm.MAIN.off}.
It is easy to see that  
$\alpha_0 = {\mathbb E}[ e^{\la \bar{\nu},  \phi_1^2\ra }  \delta (\phi_1^1)] >0$.
\toy


\begin{thebibliography}{99}

\bibitem{ai}
Aida, S.;  Vanishing of one-dimensional $L^2$-cohomologies of loop groups. 
J. Funct. Anal. 261 (2011), no. 8, 2164--2213.




\bibitem{az}
Azencott, R.;
Densit\'e des diffusions en temps petit: 
d\'eveloppements asymptotiques. I.  
Seminar on probability, XVIII, 402--498,
Lecture Notes in Math., 1059, Springer, Berlin, 1984. 


\bibitem{be}
Bailleul, I.;
Flows driven by rough paths. 
 Rev. Mat. Iberoam. 31 (2015), no. 3, 901--934. 
 
 
\bibitem{be2}
Bailleul, I.;
Regularity of the Ito-Lyons map. 
Confluentes Math. 7 (2015), no. 1, 3--11. 


\bibitem{bnot}
Baudoin, F.; Nualart, E.; Ouyang, C; Tindel, S.;
On probability laws of solutions to differential systems driven by a fractional Brownian motion.
To appear in Ann. Probab.


\bibitem{bo_svy}
Baudoin, F.; Ouyang, C.; 
On small time asymptotics for rough differential equations driven by fractional Brownian motions.
In Large deviations and asymptotic methods in finance, 413--438.  
Springer, 2015.

\bibitem{boz1}
Baudoin, F.; Ouyang, C.;  Zhang, X.;
Varadhan Estimates for rough differential equations driven by fractional Brownian motions.
Stochastic Process. Appl. 125 (2015), Issue 2, 634--652. 
%

\bibitem{boz2}
Baudoin, F.; Ouyang, C.;  Zhang, X.;
Smoothing effect of rough differential equations driven by fractional Brownian motions.
 Ann. Inst. Henri Poincar\'e Probab. Stat. 52 (2016), no. 1, 412--428.


\bibitem{bena1}
Ben Arous, G.;
Noyau de la chaleur hypoelliptique et g\'eom\'etrie sous-riemannienne. 
Stochastic analysis (Paris, 1987), 1--16,

\bibitem{bena2}
Ben Arous, G.;
D\'eveloppement asymptotique du noyau de la chaleur hypoelliptique hors du cut-locus. 
Ann. Sci. \'Ecole Norm. Sup. (4) 21 (1988), no. 3, 307--331.

\bibitem{bena3}
Ben Arous, G.;
D\'eveloppement asymptotique du noyau de la chaleur hypoelliptique sur la diagonale. 
Ann. Inst. Fourier (Grenoble) 39 (1989), no. 1, 73--99.




\bibitem{bl1}
Ben Arous, G.; L\'eandre, R.;
D\'ecroissance exponentielle du noyau de la chaleur sur la diagonale. I. 
 Probab. Theory Related Fields 90 (1991), no. 2, 175--202. 

\bibitem{bl2}
Ben Arous, G.; L\'eandre, R.;
D\'ecroissance exponentielle du noyau de la chaleur sur la diagonale. II. 
 Probab. Theory Related Fields 90 (1991), no. 3, 377--402. 





\bibitem{bi}
Bismut, J.-M.;
Large deviations and the Malliavin calculus. Progress in Mathematics, 45. Birkh\"auser Boston, Inc., Boston, MA, 1984.

\bibitem{bgq}
Boedihardjo, H.; Geng, X.; Qian, Z.;
Quasi-sure existence of Gaussian rough paths and large deviation principles for capacities.
To appear in Osaka J. Math.

 
\bibitem{cf}
 Cass, T.; Friz, P.; 
 Densities for rough differential equations under H\"ormander's condition. 
 Ann. of Math. (2) 171 (2010), no. 3, 2115--2141. 

\bibitem{cfv}
 Cass, T.; Friz, P.; Victoir, N.; 
 Non-degeneracy of Wiener functionals arising from rough differential equations. Trans. Amer. Math. Soc. 361 (2009), no. 6, 3359--3371.



\bibitem{chlt}
Cass, T.; Hairer, M.; Litterer, C.;  Tindel, S.;
Smoothness of the density for solutions to Gaussian Rough Differential Equations.
 Ann. Probab. 43 (2015), no. 1, 188--239.
 
\bibitem{cll}
Cass, T.; Litterer, C.; Lyons, T.;
 Integrability and tail estimates for Gaussian rough differential equations. 
 Ann. Probab. 41 (2013), no. 4, 3026--3050. 

\bibitem{dri}
 Driscoll, P.;
 Smoothness of densities for area-like processes of fractional Brownian motion. Probab. Theory Related Fields 155 (2013), no. 1-2, 1--34.
 
 
 \bibitem{fggr}
 Friz, P; Gess. B.; Gulisashvili, A.; Riedel, S.;
 Jain-Monrad criterion for rough paths and applications.
 Ann. Probab. 44 (2016), no. 1,  684 - 738.
    
 \bibitem{fo}
Friz, P.; Oberhauser, H;  
A generalized Fernique theorem and applications.  
Proc. Amer. Math. Soc.  138  (2010),  no. 10, 3679--3688.

 
 
 
 \bibitem{fv06}
Friz, P.; Victoir, N.;
A variation embedding theorem and applications.  J. Funct. Anal.  239  (2006),  no. 2, 631--637. 
 
 \bibitem{fv07}
Friz, P.; Victoir, N.;
Large deviation principle for enhanced Gaussian processes.
 Ann. Inst. H. Poincar\'e Probab. Statist. 43 (2007), no. 6, 775--785.  
 
 




\bibitem{fvbk}
Friz, P.; Victoir, N.;
Multidimensional stochastic processes as rough paths.
 Cambridge University Press, Cambridge, 2010. 


\bibitem{gav}
Gaveau, B.;
Principe de moindre action, propagation de la chaleur et estim\'ees 
sous elliptiques sur certains groupes nilpotents.
Acta Math. 139 (1977), no. 1-2, 95--153. 


\bibitem{hp}
 Hairer, M.; Pillai, N.; 
 Regularity of laws and ergodicity of hypoelliptic SDEs driven by rough paths. 
 Ann. Probab. 41 (2013), no. 4, 2544--2598

\bibitem{ht}
 Hu, Y.; Tindel, S.; 
 Smooth Density for Some Nilpotent Rough Differential Equations. 
 J. Theoret. Probab. 26 (2013), no. 3, 722--749.

\bibitem{iwbk}
Ikeda, N.; Watanabe, S.;
Stochastic differential equations and diffusion processes. Second edition. 
North-Holland Publishing Co., Amsterdam; Kodansha, Ltd., Tokyo, 1989.





\bibitem{ina06}
Inahama, Y.; 
Quasi-sure existence of Brownian rough paths and a construction of Brownian pants. Infin. Dimens. Anal. Quantum Probab. Relat. Top. 9 (2006), no. 4, 513--528.  




\bibitem{ina10}
Inahama, Y.;
A stochastic Taylor-like  expansion in the rough path theory. 
J. Theoret. Probab. 23 (2010) 671--714.

\bibitem{inaaop}
Inahama, Y.;
Laplace approximation for rough differential equation
driven by fractional Brownian motion, 
 Ann. Probab. 41 (2013), No. 1, 170-205.


\bibitem{ina1}
Inahama, Y.;
Short time kernel asymptotics for Young SDE 
by means of Watanabe distribution theory.
J. Math. Soc. Japan 68, No. 2 (2016), 1--43. 

\bibitem{ina2}
Inahama, Y.;
Large deviation principle of Freidlin-Wentzell type for pinned diffusion processes.
Trans. Amer. Math. Soc. 367 (2015), 8107-8137.
%
 
\bibitem{ina3}
Inahama, Y.;
Malliavin differentiability of solutions of rough differential equations.
J. Funct. Anal. 267 (2014), 1566-1584. 



\bibitem{ik}
Inahama, Y.; Kawabi, H,;
Asymptotic expansions for the Laplace approximations for It\^o functionals of Brownian rough paths.  J. Funct. Anal.  243  (2007),  no. 1, 270--322.



\bibitem{kuo}
Kuo, H.-H.;
Gaussian measures in Banach spaces.
Lecture Notes in Mathematics, Vol. 463. Springer-Verlag, Berlin-New York, 1975.





\bibitem{ks1}
Kusuoka, S.; Stroock, D. W.;
Precise asymptotics of certain Wiener functionals.  
J. Funct. Anal.  99  (1991),  no. 1, 1--74. 

\bibitem{ks2}
Kusuoka, S.; Stroock, D. W.;
Asymptotics of certain Wiener functionals with degenerate extrema.  Comm. Pure Appl. Math.  47  (1994),  no. 4, 477--501.






\bibitem{lean1}
L\'eandre, R.;
Majoration en temps petit de la densit\'e d'une diffusion d\'eg\'en\'er\'ee. 
Probab. Theory Related Fields 74 (1987), no. 2, 289--294.



\bibitem{lean2}
L\'eandre, R.;
Minoration en temps petit de la densit\'e d'une diffusion  d\'eg\'en\'er\'ee.  
J. Funct. Anal. 74 (1987), no. 2, 399--414. 

\bibitem{lean3}
L\'eandre, R.;
Int\'egration dans la fibre associ\'ee \`a  une diffusion d\'eg\'en\'er\'ee. 
Probab. Theory Related Fields 76 (1987), no. 3, 341--358.


\bibitem{lean4}
L\'eandre, R.;
Applications quantitatives et g\'eom\'etriques du calcul de Malliavin. 
Stochastic analysis (Paris, 1987), 109--133, 
Lecture Notes in Math., 1322, Springer, Berlin, 1988.


\bibitem{lean5}
L\'eandre, R.;
D\'eveloppement asymptotique de la densit\'e d'une diffusion d\'eg\'en\'er\'ee. 
Forum Math. 4 (1992), no. 1, 45--75.

\bibitem{le}
 Lejay, A.; 
 An introduction to rough paths. S\'eminaire de Probabilit\'es 
 XXXVII, 1--59, Lecture Notes in Math., 1832, Springer, Berlin, 2003.
 
\bibitem{lcl}
 Lyons, T.; Caruana, M.; L\'evy, T.; 
Differential equations driven by rough paths. 
  Lecture Notes in Math., 1908. Springer, Berlin, 2007.

\bibitem{lq}
 Lyons, T.; Qian, Z.; 
 System control and rough paths. 
 Oxford University Press, Oxford, 2002. 

\bibitem{mol}
Molchanov, S. A.;
Diffusion processes, and Riemannian geometry.
Russian Math. Surveys 30 (1975), no. 1, 1--63.


\bibitem{nu}
 Nualart, D.; 
 The Malliavin calculus and related topics. 
 Second edition. Springer-Verlag, Berlin, 2006. 



\bibitem{sh}
Shigekawa, I.;
Stochastic analysis.
Translations of Mathematical Monographs, 224.
 Iwanami Series in Modern Mathematics. American Mathematical Society, Providence, RI, 2004.


\bibitem{tak}
 Takanobu, S.; 
 Diagonal short time asymptotics of heat kernels for certain degenerate
  second order differential operators of H\"ormander type. 
 Publ. Res. Inst. Math. Sci. 24 (1988), no. 2, 169--203. 

\bibitem{tw}
Takanobu, S.; Watanabe, S.;
Asymptotic expansion formulas of the Schilder type for a class of conditional Wiener functional integrations. 
Asymptotic problems in probability theory: Wiener functionals and asymptotics (Sanda/Kyoto, 1990), 
194--241, Pitman Res. Notes Math. Ser., 284, Longman Sci. Tech., Harlow, 1993.


\bibitem{ue1}
Uemura, H.; 
On a short time expansion of the fundamental solution of heat equations by the method of Wiener functionals. J. Math. Kyoto Univ. 27 (1987), no. 3, 417--431. 

\bibitem{ue2}
Uemura, H.; Off-diagonal short time expansion of the heat kernel on a certain nilpotent Lie group.
J. Math. Kyoto Univ. 30 (1990), no. 3, 403--449.


\bibitem{uw}
Uemura, H.; Watanabe, S.;
Diffusion processes and heat kernels on certain nilpotent groups. 
Stochastic analysis (Paris, 1987), 173--197, Lecture Notes in Math., 1322, Springer, Berlin, 1988.

\bibitem{wa}
Watanabe, S.; 
Analysis of Wiener functionals (Malliavin calculus) and its applications to heat kernels.  
Ann. Probab.  15  (1987),  no. 1, 1--39.




\end{thebibliography}
\end{document}